\documentclass[11pt]{article}
\usepackage[left=2.5cm, right=2.5cm, top= 1.7cm, bottom=2cm]{geometry}
\usepackage[utf8]{inputenc}
\usepackage{aecompl}
\usepackage{ae}
\usepackage{textcomp}
\usepackage[english]{babel}
\usepackage{graphicx}
\usepackage{pdfsync}
\usepackage{tikz}
\usetikzlibrary{arrows,positioning}
\usetikzlibrary{calc}
\usepackage{pgffor}
\usepackage{placeins}
\usepackage{scalefnt}
\usepackage{subfigure}
\usepackage{caption}
\usepackage{hyperref}
\hypersetup{
    colorlinks,%
    citecolor=blue,%
    filecolor=black,%
    linkcolor=blue,%
    urlcolor=black
}
\usepackage{hypbmsec}
\usepackage{amsmath,amsthm,amsfonts,amssymb}
\usepackage{bbm}
\usepackage{enumerate}
\usepackage{dsfont}
\usepackage{url}
\usepackage{float}
\usepackage{setspace}

\def\R{\mathbb{R}}

\def\Z{\mathbb{Z}}
\def\G{\mathsf{G}}
\def\C{\mathcal{C}_2}
\def\B{\mathcal{B}}
\newcommand{\SG}{\mathsf{SG}}
\newcommand{\SC}{\mathsf{SC}}

\newcommand{\abs}[1]{\lvert#1\rvert}

\theoremstyle{plain}
\newtheorem{theorem}{Theorem}[section]

\newtheorem{definition}[theorem]{Definition}
\newtheorem{question}[theorem]{Question}
\newtheorem{conjecture}[theorem]{Conjecture}

\theoremstyle{definition}

\title{\bf From fractals in external DLA to internal DLA on fractals}
\author{Ecaterina Sava-Huss \footnote{Institute of Discrete Mathematics, Graz University of Technology, Austria. sava-huss@tugraz.at; http://www.math.tugraz.at/~sava}}

\begin{document}
\maketitle

\begin{abstract}
We present an unified approach on the behavior of two random growth models (\textit{external DLA} and \textit{internal DLA}) on infinite graphs, the second one being an internal counterpart of the first one. Even though the two models look pretty similar, their behavior is completely different: while external DLA tends to build irregularities and fractal-like structures, internal DLA tends to fill up gaps and to produce regular clusters. We will also consider the aforementioned models on fractal graphs like Sierpinski gasket and carpet, and present some recent results and possible questions to investigate. 
\end{abstract}

\begin{spacing}{1.0}
\begin{small}
\tableofcontents
\end{small}
\end{spacing}

\setlength{\parskip}{2.0ex}

\section{Introduction}
We consider two aggregation models initially introduced in physics in  \cite{dla} and \cite{meakin-deutch},
and rigorously studied in mathematics over the last three decades, models for which we present a survey on the existing results and state several open problems. The models under consideration are \textit{external diffusion limited aggregation} (shortly \textit{external DLA}) and \textit{internal diffusion limited aggregation} (shortly \textit{internal DLA}). In the mathematical community, these two models started to gain interest only  a couple of years after being introduced, with the first results on external DLA in \cite{kesten-dla-hit-pro,kesten-dla}, and on internal DLA in \cite{lawler_bramson_griffeath}. Only recently, these models became interesting in the fractals community:
few recent results concerning external DLA on the $m$-dimensional pre-Sierpinski carpet as defined in \cite{osada-pre-sierpinski}, for $m\geq 3$ are available. For the internal DLA
on the Sierpinski gasket graph, there are  also some limit shape results, but other than these two examples, there is not much known about the two growth models on other fractal graphs, where according to simulations which we present towards to end of the paper, interesting behavior may be observed. With the current overview, we would like to draw  the attention on the beauty of these models.

For the rest of the paper, $\G$ will be an infinite and locally finite graph, the reference state space, which will be replaced with concrete examples of graphs as needed. We denote by $o\in \G$
a fixed vertex, the \emph{origin of the graph} $\G$.

\textbf{External DLA} was initially introduced in physics by {\sc Witten and Sander} \cite{dla} as an example to create ordering out of chaos due to a simple rule. Mathematically, this ordering is far away from being understood, and new methods and ideas are needed in order to move forward in this direction. External DLA
 is a model of random fractal growth which exhibits self-organized criticality and complex-pattern formation, and which produces scale-invariant objects whose Hausdorff dimension is independent of short-range details. Moreover external DLA has no upper critical dimension as shown in \cite{dla}; it is a model which builds a sequence of random growing sets $(\mathcal{E}_n)_{n\geq 0}$, starting with one particle $\mathcal{E}_0=\{o\}$ at the origin of $\G$. At each time step, a new particle starts a simple random walk from "infinity" (far away) and walks until it hits the outer boundary of the existing cluster, where it stops and settles. In this way, one builds a family $(\mathcal{E}_n)_{n\geq 0}$ of growing clusters; the set $\mathcal{E}_n$ consists of exactly $n+1$ particles and it is called \emph{external DLA cluster}. In spite of these very simple growth rules, only a few rigorous mathematical results about
external DLA are available, results which will be surveyed below. A typical structure produced on a two-dimensional lattice is shown in Figure \ref{fig:dla}. External DLA was found to well represent growth processes in nature such as growth of bacterial colonies, electrodeposition, or crystal growth.
\begin{figure}
\label{fig:dla}
\begin{center}
\includegraphics[height=6cm]{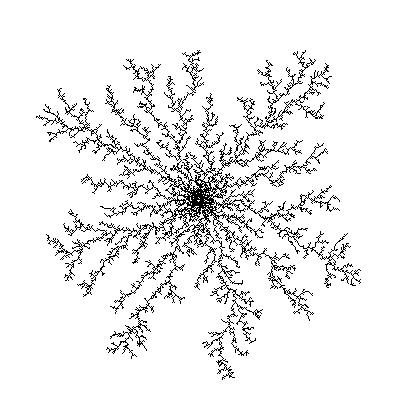}
\end{center}
\caption{\label{fig:dla}: External DLA cluster on $\mathbb{Z}^2$ with center initially occupied}
\end{figure}

\textbf{Internal DLA} is an attempt of a model which eliminates irregularities and fills gaps, as opposed to external DLA. It was proposed by {\sc Meakin and Deutch} \cite{meakin-deutch} as a model of industrial chemical processes such as electropolishing, corrosion and etching. {\sc Diaconis and Fulton} in \cite{diaconis_fulton_1991} identified internal DLA as a special case of a “smash sum” operation on subsets of $\Z^2$. Internal DLA is a random growth model which builds a sequence of random growing clusters $(\mathcal{I}_n)_{n\geq 0}$ based on particles performing random walks, where all the particles start from the same fixed point $o$. Typically, one starts with $\mathcal{I}_0=\{o\}$, and for each $n$, we let $\mathcal{I}_{n+1}$ be $\mathcal{I}_n$ plus the first point where a random walk started at $o$ exits $\mathcal{I}_n$. There are several modifications in this model, where one can start the random walks uniformly at random in the already existing cluster, or one can start with an initial configuration of particles on the state space $\G$. Results in this directions will be surveyed in the following. As in external DLA, understanding the shape of the limiting cluster $\mathcal{I}_n$, the \textit{internal DLA cluster} with $n+1$ particles, is the main question in this model.
Also, of fundamental significance as mentioned in the initial paper \cite{meakin-deutch}, is to know how smooth a surface formed by internal DLA (processes) may be. These problems are well understood mathematically on many state spaces, and there are very precise results. 
On the one hand, the limiting object formed from internal DLA does not show any fractal structure. On the other hand, when running internal DLA on a fractal graph, we have partial results that indicate the absence of fractal structure, though there remain many more fractal state spaces to be explored.

The crucial difference between the above two models is that the dynamics of the external model roughens the cluster, whereas the dynamics of the internal model makes the cluster smoother.

\textbf{Structure of the paper.} After fixing the notation and the basic notions in Section \ref{sec:prel}, we focus on the external DLA model in Section \ref{sec:edla}, in which we survey the available results on the growth of arms in this model, number of holes, and variations of the standard model. The results will not be stated in the chronological order of publication, but according to the state space they evolve on. Finally, in Section \ref{sec:idla} we survey the results for the limit shapes of the internal DLA cluster, and we include several questions through the whole paper. 

\section{Preliminaries}\label{sec:prel}

\paragraph{Graphs.} Let $\G$ be an infinite, locally finite graph (i.e. every vertex has finite degree denoted by $\deg(x)$). The neighborhood relation will be denoted by $"\sim"$, and by $x\sim y$ we mean that $(x,y)$ is an edge in $\G$. Let $o\in \G$ be a fixed distinguished vertex, which will be called the \emph{origin} or the \emph{root}. For $x,y\in \G$, the distance $d(x,y)$ represents the minimal number of edges on the path connecting $x$ with $y$. For a subgraph $A$ of $\G$, we denote by $\partial A$ the outer boundary of $A$:
$$\partial A=\{y\in \G: y\notin A,\ \exists x\in A: \ x\sim y\}.$$
For $x\in \G$ and $n\geq 0$, we write $\mathcal{B}_n(x)=\{y\in G:\ d(x,y)\leq n\}$ for the ball of radius $n$ and center $x$ in $\G$. If the center of the ball is $o$, we write only $\mathcal{B}_n$.

\paragraph{Random Walks.} Let $(S_n)_{n\geq 0}$ be a random walk on $\G$, and denote by $\mathbb{P}_x$ the probability measure of the random walk started at $x$. We do not fix yet the transition probabilities for the random walk, since those will change from case to case, and we will mention them as needed. For a subset $A\subset \G$, let $T(A)$ be the hitting time of $A$, defined as
$$T(A)=\min\{n\geq 0:\ S_n\in A\}.$$
For a set $\{x\}$ with a single vertex, we write $T(x)$ instead of $T(\{x\})$. The \emph{heat kernel} of the random walk $S_n$ is defined to be
$$p_n(x,y)=\mathbb{P}_x[S_n=y],$$ and the \emph{Green function} $G(x,y)$ is defined as
$$G(x,y)=\sum_{n\geq 0}p_n(x,y),$$
which is well defined and finite precisely when the random walk is transient. For a subset $A\subset \G$, the \emph{killed or the stopped Green function} $G_A(x,y)$ is defined as
$$G_A(x,y)=\sum_{n\geq 0}\mathbb{P}_x[S_n=y,T(A)>n].$$ The \emph{hitting distribution} $H_A(x,y)$ is then
$$H_A(x,y)=\mathbb{P}_x[S_{T(A)}=y],\quad \text{for }y\in A$$
and $S_{T(A)}$ is the hitting position of $A$. If the random walk $(S_n)$ starts at $o$,
we write
\begin{equation}\label{eq:harm-meas}
h_{A}(y)=\mathbb{P}_o[S_{T(A)}=y],\quad \text{for } y\in A
\end{equation}
for the probability of the random walk starting at $o$ to first hit $A$ in $y$, that is $h_A$ is the harmonic measure (from $o$) of the set $A$, and $\sum_{y\in A}h_A(y)=1$.

\section{External DLA}\label{sec:edla}

We define here formally the external DLA model, by first explaining what it means to release a particle at infinity. Several variants of external DLA have been considered, but we refer here to the original, simplest model, which can be defined on any space
where the notion of random walk or diffusion  exists.  If 
the Poisson boundary consists of one point  and 
 the random walk is recurrent (for instance the case of simple random walks on $\Z$ and $\Z^2$), external DLA  can be defined so that the law of the location of a new particle is the harmonic
measure of the existing aggregate with pole at infinity. If the random walk is
transient (such as the case of simple random walks on $\Z^d$, with $d\geq 3$, or on regular trees $\mathbb{T}_d$ of degree $d\geq 2$, $n$-dimensional Sierpinski carpet graph, for $n\geq 3$), one can consider the harmonic measure
with a pole far away from the aggregate, let the pole go to infinity and take
limits (i.e., conditioning the random walk coming from infinity to hit the
cluster). 
That is, in defining rigorously external DLA, we have to distinguish the cases when the random walk $(S_n)$ on the infinite graph $\G$ is recurrent or transient; the Poisson boundary of the random walk also plays a role in this case.
We recall that the {\it Poisson boundary} of a random walk is measure space that describes the stochastically significant behavior of the walk at infinity. It provides an integral representation of the bounded harmonic functions of the random walk.

During the whole paper, when we speak about the $n$-dimensional Sierpinski carpet graph, we shall also use the notion {\it pre-Sierpinski carpet}, and we have in mind the construction introduced in  \cite{osada-pre-sierpinski}.

We shall write $\mu_A(y)$ for the \emph{harmonic measure from infinity}, that is, for the probability to start a random walk at infinity and to hit the finite subset $A\subset \G$ at the point $y$. Depending on whether the graph $\G$ is transient or recurrent, this measure can take different forms, and we cannot define it globally on any general graph here. This will be made precise in the concrete cases below.

\begin{definition}\label{def:dla-mc} Let $\G$ be an infinite graph, and $(S_n)$ a random walk on it. \emph{External DLA} on $\G$ is a Markov chain $(\mathcal{E}_n)_{n\geq 0}$ on finite subsets of $\G$, which evolves in the following way. Start with a single vertex $o\in \G$, that is $\mathcal{E}_0=\{o\}$. Given the process $\mathcal{E}_n$
at time $n$, let $y_{n+1}$ be a random vertex in $\partial \mathcal{E}_n$ chosen according to the harmonic measure (from infinity) of $\partial \mathcal{E}_n$.
That is,
$$\mathbb{P}[y_{n+1}=y|\mathcal{E}_n]=\mu_{\partial \mathcal{E}_n}(y),\quad \text{for }y\in\partial \mathcal{E}_n,$$ 
and we set $\mathcal{E}_{n+1}=\mathcal{E}_{n}\cup \{y_{n+1}\}.$ 
\end{definition}

\begin{definition}
The \emph{cluster at infinity $\mathcal{E}_{\infty}$} for the external DLA process $(\mathcal{E}_n)$ on $\G$ is defined as
$$\mathcal{E}_{\infty}=\bigcup_{n=1}^{\infty}\mathcal{E}_{n}.$$
\end{definition}

It is immediate that the external DLA cluster $\mathcal{E}_n$ at time $n$ contains exactly $n+1$ particles.
This model is hard to study. The difficulty comes from the fact that
the dynamics is neither monotone nor local (meaning that if big tentacles surround a vertex $x$, than $x$ will never be added to the cluster). By non-monotonicity we mean that there
is no coupling between the external DLA starting from a cluster $C$ and another from a
cluster $D\subset C$ such that, at each step, the inclusion of the clusters remains valid almost surely.
Understanding the shape of $\mathcal{E}_n$ as $n\to \infty$ and the fractal nature of this object, are problems one would be typically interested in. While mathematically this is out of reach for the time being, there are other partial results concerning the growth of arms and the number of holes in external DLA.

\subsection{Integer lattices $\Z^d$}

In this subsection the state space for the external DLA process is $\G=\Z^d$, $d\geq 1$. Even for $\Z^d$, there are no results that prove the fractal nature of the limiting object, or results that prove the zero density in the long run. The first rigorous results go back to Kesten \cite{kesten-dla-hit-pro,kesten-dla}, who gives estimates on the growth of arms in external DLA.
Since for $d=1$, the behavior of standard external DLA is trivial, we consider $d\geq 2$, and let $(S_n)_{n\geq 0}$ be a simple random walk on $\Z^d$.

\textbf{For $d=2$,} for any finite nonempty subset $A\subset \Z^2$, we have $T(A)<\infty$ with probability one, and we define the \emph{harmonic measure (from infinity) of $A$}
\begin{equation}\label{eq:harm-measure-rec}
\mu_A(y)=\lim_{|x|\to\infty} H_A(x,y),
\end{equation}
where $|x|$ denotes the Euclidean norm of $x$. The limit $\lim_{|x|\to\infty}$ corresponds to "releasing the particle at infinity". In this case, $(S_n)$ is recurrent, so that by \cite[Theorem 14.1]{spitzer-principles} the limit in \eqref{eq:harm-measure-rec} exists and $\sum_{y\in A}\mu_A(y)=1$.

\textbf{For $d\geq 3$}, since the random walk $(S_n)$ is transient, the limit $\lim_{|x|\to\infty} H_A(x,y)$ in \eqref{eq:harm-measure-rec} is identically zero (cf \cite[Proposition 25.3]{spitzer-principles}). So in order to obtain a nontrivial limit similar to the one in \eqref{eq:harm-measure-rec}, we have to condition on $T(A)$ being finite. This conditioning
gives the factor of the capacity of the set $A$ in the denominator.
In the case $d\geq 3$, we define the \emph{harmonic measure (from infinity)} of a finite subset $A\subset G$ as
\begin{equation}\label{eq:harm-measure-trans}
\mu_A(y)=\lim_{d(o,x)\to\infty}\dfrac{H_A(x,y)}{\sum_{z\in A}H_A(x,z)}=\lim_{d(o,x)\to\infty}\mathbb{P}_x[S_{T(A)}=y|T(A)<\infty],\quad \text{for }y\in A,
\end{equation}
which is proportional to the so-called \emph{equilibrium measure} associated to the
set $A$. The limit in \eqref{eq:harm-measure-trans} exists again by \cite[Proposition 26.2]{spitzer-principles} for $d=3$ (the same proof works also for $d>3$) and satisfies $\sum_{y\in A}\mu_A(y)=1$. Therefore, we have a valid definition for external DLA, and we let $r(\mathcal{E}_n)$ to be the radius of $\mathcal{E}_n$, defined as
\begin{equation}\label{eq:radius-edla}
r(\mathcal{E}_n)=\max\{|x|:\ x\in \mathcal{E}_n\},
\end{equation}

\begin{theorem}{\cite[Theorem]{Kesten-phys} and \cite[Corollary]{kesten-dla-hit-pro}}
\label{thm:kesten-dla}
There exist constants $C(d)<\infty$ such that with probability $1$
\begin{align*}
\limsup_{n\to\infty}n^{-2/3}r(\mathcal{E}_n) & \leq C(2),\quad \text{ if }d=2\\
\limsup_{n\to\infty}n^{-2/d}r(\mathcal{E}_n) & \leq C(d),\quad \text{ if }d\geq 3.
\end{align*}
\end{theorem}
The proof uses classical estimates for the harmonic measure (from infinity) as defined in  \eqref{eq:harm-measure-rec} and \eqref{eq:harm-measure-trans} and for the hitting probabilities.
Simulations actually indicate that for $d=2$, $\mathbb{E}[r(\mathcal{E}_n)]\approx n^{10/17}$ but 
 as far as the lower bound is concerned, nothing has been proven beyond $\sqrt{n}$ in the 35 years since the model has been introduced. It would be very
interesting to prove even a logarithmic correction, i.e. to prove that $\mathbb{E}[r(\mathcal{E}_n)]\geq \sqrt{n}\log(n)$.
On $\Z^d$, a lower bound on the number $N(n)$ of vertices in $\B_n$ which are occupied by the cluster $E_{\infty}$ is known.
\begin{theorem}{\cite[Theorem 2]{kesten-dla}}
There exist constants $C(d)$ such that with probability $1$
$$N(n)\geq C(d)n^{d-1},\quad \text{ for infinitely many }n.$$
\end{theorem} 
Another non-trivial result on $\Z^2$ concerns the number of holes $\mathcal{H}_n$ in the external DLA cluster $\mathcal{E}_n$. A hole of $\mathcal{E}_n$ is a finite connected component of $\Z^2\setminus \mathcal{E}_n$. 
\begin{theorem}{\cite{eberz-wagner-dla}} For any finite connected subset $e$ of $\Z^2$ we have
$$\mathbb{P}[\mathcal{H}_n \text{ converges to infinity as } n \to\infty | \mathcal{E}_0 = e] = 1.$$
\end{theorem}
Theorem \ref{thm:kesten-dla} has been improved in \cite{benjamini-yadin-dla}, where upper bounds on the growth rate of arms in external DLA cluster are given on a big class of transient graphs such as: transitive graphs of polynomial growth of degree $\geq 4$; transitive graphs of exponential growth; $\mathbb{Z}^3$; non-amenable graphs; $n$-dimensional pre-Sierpinski gasket graphs ($n\geq 3$) as introduced in \cite{osada-pre-sierpinski}. In particular, on $\Z^3$ the factor $n^{-2/3}$ from Theorem \ref{thm:kesten-dla} has been improved to $n^{-1/2}/\log(n).$ On the class of transient graphs $\G$ considered in \cite{benjamini-yadin-dla}, the harmonic measure (from infinity) $\mu_A$ of a set $A$ is defined as in \eqref{eq:harm-measure-trans}.

A directed version of external DLA has been recently introduced on $\Z^2$ in \cite{martineau-directed-dla}. In a series of three papers \cite{dla-long-range1,dla-long-range2, amir-dla}, a one-dimensional external DLA model based on random walks with long jumps (that depend on a parameter $\alpha$) is proposed, which tries to capture the fractal nature of the standard DLA. Depending on the values of $\alpha$, the random walk $(S_n)$ with long jumps on $\Z$ may be recurrent or transient, and for the precise definition of harmonic measure from infinity we refer to those three papers.
The main results of \cite{dla-long-range1,dla-long-range2} can be summarized into the following theorem.
\begin{theorem}
Let $(S_n)$ be a symmetric random walk on $\Z$ that satisfies
$\mathbb{P}[|S_1-S_0|=k]\sim ck^{-1-\alpha}$. Let $d(\mathcal{E}_n)$ be the diameter of the external DLA cluster $\mathcal{E}_n$. Then almost surely:
\hspace{-5em}
\begin{enumerate}[(a)]
\setlength\itemsep{0em}
\item If $\alpha>3$, then $n-1\leq d(\mathcal{E}_n)\leq Cn+o(n)$, where $C$ depends only on $\alpha$.
\item If $2<\alpha\leq 3$, then $d(\mathcal{E}_n)=n^{\beta+o(1)}$, where $\beta=\frac{2}{\alpha-1}$.
\item If $1<\alpha<2$, then $d(\mathcal{E}_n)=n^{2+o(1)}$.
\item If $\frac{1}{3}<\alpha<1$, then $n^{\beta+o(1)}\leq d(\mathcal{E}_n)\leq n^{\beta'+o(1)}$, where $\beta=\max(2,\alpha^{-1})$ and $\beta'=\frac{2}{\alpha(2-\alpha)}$.
\item If $0<\alpha<\frac{1}{3}$, then $d(\mathcal{E}_n)=n^{\beta+o(1)}$, where $\beta=\alpha^{-1}$.
\end{enumerate}
\end{theorem}
The last one \cite{amir-dla} from the series of three papers mentioned above deals with the cluster at infinity $E_{\infty}$, and it is shown that for random walks $(S_n)$ whose step size has finite third moment, $E_{\infty}$ has a
renewal structure and positive density. In contrast, for random walks whose step size has finite variance, the renewal structure no longer exists and $E_{\infty}$ has zero density. 

\begin{theorem}{\cite[Theorem 1]{amir-dla}} Assume that the step distribution $\xi$ of the random walk $(S_n)$ on $\Z$ satisfies $\mathbb{P}[\xi>n]\leq Cn^{-\alpha}$ for any $n$ and some $\alpha>3$. There exists some $B>0$ such that a.s. $\mathcal{E}_{\infty}$ has density $B$. Further, $B$ is the limit density of $\mathcal{E}_n$:
$$B=\lim_{m_1\to\infty,m_2\to\infty}\dfrac{|\mathcal{E}_{\infty}\cap[-m_1,m_2]|}{m_1+m_2}=\lim_{n\to\infty}\frac{n}{d(\mathcal{E}_n)}.$$
\end{theorem}
\begin{theorem}{\cite[Theorem 2]{amir-dla}}
Assume that there exist $2<\alpha<3$ and constants $c_1,c_2>0$ so that $\xi$ satisfies $c_1n^{-\alpha}\leq \mathbb{P}[\xi>n]\leq c_2n^{-\alpha}$ for all $n$ then a.s.
$$|\mathcal{E}_{\infty}\cap[-n,n]|=n^{\frac{\alpha-1}{2}+o(1)}.$$
In particular, $\mathcal{E}_{\infty}$ has zero density in the sense that $\lim_{n\to\infty}\frac{|\mathcal{E}_{\infty}\cap[-n,n]|}{n}=0$.
\end{theorem}
The results mentioned above are the only ones available for external DLA on $\Z^d$, and the limit shape and the density problem for $d\geq 2$ still resist a mathematical proof. There are many open problems and questions in this direction; see \cite{benjamini-yadin-dla} for more details.

\begin{conjecture} On $\Z^d$, the rate of growth of the radius of the external DLA cluster $\mathcal{E}_n$
started at $\mathcal{E}_0=\{0\}$ is of order $n^{1/d}$:
$$\limsup_{n\to\infty}n^{-1/d}\mathbb{E}[r(\mathcal{E}_n)]=0.$$
\end{conjecture}
\begin{question}
What is the distribution of the number of ends of the cluster at infinity $\mathcal{E}_{\infty}$ on $\mathbb{Z}^d$?
\end{question}
Concerning recent progress on external DLA in a wedge of $\mathbb{Z}^d$, we refer to \cite{dla-wedge}. Furthermore, the reach of Kesten’s idea is extended to non-transitive
graphs in \cite{Procaccia2019}, where the (horizontally) translation invariant stationary harmonic
measure on the upper half plane with absorbing boundary condition is defined and it is shown that the
growth of such stationary harmonic measure in a connected subset intersecting
x-axis is sub-linear with respect to the height; see also \cite{procaccia-zero-harmonic,procaccia-scaling-limit} where the stationary harmonic measure as a natural growth measure for external DLA model in the upper planar lattice is investigated.

\subsection{Trees $\mathbb{T}_d$}

One reason that makes the lattice case $\Z^d$ hard to investigate is that there is no simple way to describe the harmonic measure (from infinity) for the boundary of an external DLA cluster on $\Z^d$. On other state spaces, such as trees, which have no loops, the model is more tractable and the harmonic measure (from infinity) can be understood. In \cite{barlow-pemantle-perkins-dla}, an adjusted version of external DLA on $d$-regular trees $\mathbb{T}_d$, where the fingering phenomenon occurs, was introduced. The dynamics of their model is as follows:
the initial cluster $\mathcal{E}_0$ contains only the root. Vertices are then added one by one from among
those neighboring the current subtree. The choice of which vertices to add is random, with
vertices in generation $n$ (i.e. distance $n$ from the root) chosen with probabilities proportional to $\alpha^{-n}$ where $\alpha>0$ is a fixed parameter. Then $\mathcal{E}_n$ is the subtree at step $n$ and let $r(\mathcal{E}_n) = \max\{d(o,x):\ x\in \mathcal{E}_n\}$ denote the maximum height of a vertex in $\mathcal{E}_n$, which is similar to the radius in \eqref{eq:radius-edla}. For this model, for a finite subtree $A\subset \mathbb{T}_d$ with boundary $\partial A$, its harmonic measure $\mu_{\partial A}^{\alpha}$ (from infinity) on $\partial A$, with parameter $\alpha>0$ can be computed as
$$\mu_{\partial A}^{\alpha}(y)=\dfrac{\alpha^{-d(o,y)}}{\sum_{x\in\partial A}\alpha^{-d(o,y)}},\quad \text{ for } y\in\partial A,$$
see Definition on page 4 in \cite{barlow-pemantle-perkins-dla}. In the latter paper, the case $\alpha<1$ is studied. The external DLA cluster $\mathcal{E}_n$ is the position at time $n$ of the Markov chain defined in Definition \ref{def:dla-mc}, where $\mathcal{E}_{n+1}$ is obtained from $\mathcal{E}_n$ by adding a new vertex according to the harmonic measure defined in the previous equation.
For $\alpha\geq 1$ it is easy to see that $\mathcal{E}_{\infty}$ is almost surely the entire tree.
For $\alpha=1$, one has the uniform measure $\mu_{\partial A}^1$ on $\partial \mathcal{E}_n$ (this corresponds to the Eden model). From the external DLA perspective, the case $\alpha<1$ is the interesting one, where one obtains the so-called fingering phenomenon. For this external DLA model, \cite{barlow-pemantle-perkins-dla} obtained a strong law and a central limit theorem for the height $r(\mathcal{E}_n)$ of the DLA cluster.
\begin{theorem}
Let $\mathbb{T}_d$ be a d-regular tree and $0<\alpha<1$. There exist constants $r_0(\alpha,d)\in(0,1)$ and $\sigma^2=\sigma^2(\alpha,d)>0$ such that
\begin{enumerate}[(a)]
\item $\lim_{n\to\infty}\frac{r(\mathcal{E}_n)}{n}=r_0(\alpha,d)$ a.s.
\item $\frac{r(\mathcal{E}_n)-nr_0(\alpha,d)}{\sqrt{n}}\xrightarrow{\mathcal{D}} N(0,\sigma^2)$ as $n\to \infty$.
\end{enumerate}
\end{theorem}
The model considered here can be also interpreted as a model of first passage percolation on $\mathbb{T}_d$.

\subsection{Hyperbolic plane $\mathbb{H}^2$}

In \cite{eldan-dla-hyper}, external DLA on the hyperbolic plane $\mathbb{H}^2$ is considered, and it is shown 
that the cluster at infinity $E_{\infty}$ almost surely
admits a positive upper density. 
For completeness, we recall the definition of the {\it upper density} of a set, as used in \cite{eldan-dla-hyper}. In a metric measure space $X$ whose diameter is infinite, we say that a locally finite set $A\subset X$ has an {\it upper density} greater or equal to $c$ if there exist a point $p\in X$ and a sequence $R_1<R_2<\cdots$ such that $R_i\to\infty$ as $i\to\infty$, such that
$$\#\left(A\cap \mathcal{B}_{R_i}(p)\right)\geq c \mu(\mathcal{B}_{R_i}(p)),\quad \forall i\in\mathbb{N},$$
where $\mathcal{B}_{r}(p)$ is a metric ball centered at $p$ with radius $r$ and $\mu$ is the measure defined on $X$. On the hyperbolic plane, one can use this definition with the standard hyperbolic distance as a metric and the standard Riemannian volume of a set as a measure.

In the hyperbolic
setting the behavior of the aggregate is simpler to analyze than the Euclidean
one; the rate of decay of the hyperbolic potential plays an important role in understanding the external DLA. 

See Figure \ref{fig:dla-hyp} for a picture of external DLA model with 1000 particles, viewed on the Poincaré disc model.
\begin{figure}
\label{fig:dla-hyp}
\begin{center}
\includegraphics[height=4.5cm]{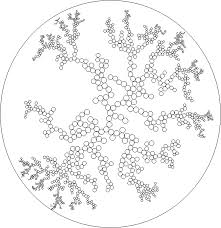}
\end{center}
\caption{\label{fig:dla-hyp}: (by Ronen Eldan) external DLA with 1000 particles, viewed on the Poincaré disc model.}
\end{figure}
In his construction, particles are metric balls of radius $1$, $\mathcal{E}_0=\{p_0\}$, where $p_0$ is a fixed point in $\mathbb{H}^2$, and recursively $\mathcal{E}_{n+1}=\mathcal{E}_{n}\cup\{y_{n+1}\}$, where $\{y_{n+1}\}$ is added to the aggregate $\mathcal{E}_n$ according to a (harmonic) measure $\mu_{\partial \mathcal{E}_n}(y)$ with pole at infinity that has to be carefully constructed on $\mathbb{H}^2$, such that external DLA makes sense in this setting. For details on this construction, we refer to \cite{eldan-dla-hyper}; the main result of his paper reads as following.
\begin{theorem}{\cite[Theorem 1.1]{eldan-dla-hyper}}
The external DLA cluster at infinity $E_{\infty}$ almost surely has an upper density greater than $c$, where $c>0$ is an universal constant. 
\end{theorem}
We would like to point out the fact that the behavior of external DLA on the hyperbolic plane and on the regular tree $\mathbb{T}_d$ as considered in \cite{barlow-pemantle-perkins-dla} is completely different, even though the hyperbolic plane has a tree-like structure.

\subsection{Cylinder graphs}
\label{sec:cyl-gr-dla}

Other results on external DLA that are worth mentioning have been proven in \cite{bejamini-yadin-cylinder} on cylinder graphs $\G\times \mathbb{N}$. Let us first fix the notation for the graphs we consider below. Let $\G$ be a finite, connected graph.
The {\it cylinder graph} with base $\G$, denoted by $\G\times \mathbb{N}$, is defined as: the vertex set of  $\G\times \mathbb{N}$ is $V(\G)\times \mathbb{N}=\{(v,k):\ v\in V(\G),\ k\in \mathbb{N}\}$, where $V(\G)$ represents the vertex set of $\G$. The edge set is defined by the following relations: for all $u,v \in V( \G)$ and all $m,k\in \mathbb{N}$,  $(u,m)\sim(v,k)$, that is between vertices $(u,m)$ and $(v,k)$ there is an edge in $\G\times \mathbb{N}$, if and only if $m=k$ and $u\sim v$ in $\G$, or $|m-k|= 1$ and  $u=v$. Equivalently, the cylinder with base $\G$ is obtained by just placing infinitely many copies of $\G$ one over the other, and connecting each vertex in a copy to its corresponding vertices in the adjacent copies.
 
On $\G\times \mathbb{N}$, particles perform simple random walks $(S_n)$ from infinity. Since $\G$ is finite, such random walks are recurrent on $\G\times \mathbb{N}$, and the harmonic measure from infinity can be defined similar to the one on $\Z^2$, as in \eqref{eq:harm-measure-rec}. That is, vertices are added to the existing cluster $\mathcal{E}_n$ according to the measure in \eqref{eq:harm-measure-rec}.
Denote by $\G_i$ the induced subgraph on the vertices of $\G\times \{m\}$, for all $m\in\mathbb{N}$, and call $\G_m$ the $m$-th level of the cylinder graph $\G\times \mathbb{N}$. One of the results proven in \cite{bejamini-yadin-cylinder} is that external DLA on $\G\times \mathbb{N}$ grows arms if the base graph $\G$ mixes fast. Recall that the mixing time $t_{\text{mix}}(\G)$ of the simple random walk on $\G$ is the time it takes for the random walk to come close in total-variation distance to the stationary distribution.
\begin{theorem}{\cite[Theorem 2.1]{bejamini-yadin-cylinder}}
Let $2\leq d\in\mathbb{N}$. There exists $n_0=n_0(d)$, such that the following holds for all $n>n_0$: let $\G$ be a $d$-regular graph of size $n$, and mixing time 
$t_{\text{mix}}(\G)\leq \dfrac{\log^2n}{(\log\log n)^5}$.
Let $(\mathcal{E}_t)$ be the external DLA process on $\G\times \mathbb{N}$ with $\mathcal{E}_0=\G_0$, and for $m\in\mathbb{N}$, let $T_m$ be the first time the DLA cluster reaches $\G_m$. Then, for all $m$, $\mathbb{E}[T_m]<\dfrac{4mn}{\log\log n}$.
\end{theorem}
This phenomenon is often referred to as \emph{the aggregate grows arms}, i.e. grows faster than order $|\G|$ particles per layer.
As mentioned in \cite{bejamini-yadin-cylinder}, the result above is believed not to be optimal, and a stronger result is conjectured.
\begin{conjecture}{\cite[Conjecture 2.2]{bejamini-yadin-cylinder}}
Let $(\mathcal{G}^n)_{n\geq 0}$ be a family of $d$-regular graphs such that $\lim_{n\to\infty}|\mathcal{G}^n|=\infty$. There exists $0<\gamma<1$ and $n_0$ such that for all $n>n_0$ the following holds: consider the cylinder graph $\mathcal{G}^n\times \mathbb{N}$ with base $\mathcal{G}^n$  and let $(\mathcal{E}_t)$ be the external DLA process on $\mathcal{G}^n\times \mathbb{N}$ with $\mathcal{E}_0$ being the zero layer of the cylinder graph, and $T_m$ be the first time the external DLA cluster reaches level $m$ on the cylinder graph $\mathcal{G}^n\times \mathbb{N}$. Then, for all $m$,
$\mathbb{E}[T_m]\leq m|\mathcal{G}^n|^{\gamma}$.
\end{conjecture} 
Concerning the density of the limit cluster at infinity $\mathcal{E}_{\infty}$, for cylinder graphs $\G\times\mathbb{N}$ with base $\G$, in the same paper there are two results. To state them, let us define the {\it empirical density of  particles in the finite cylinder} $\G\times\{1,\ldots,m\}$ as
$$D(m)=\frac{1}{mn}\sum_{i=1}^{m}|\mathcal{E}_{\infty}\cap \G_i|$$
and \emph{the density at infinity} as
$D=D_{\infty}=\lim_{m\to\infty}D(m)$. Using standard arguments from ergodic theory one can show that the above limit exists, and is constant almost surely.
 The next result relates the density at infinity to the average growth rate.
\begin{theorem}{\cite[Theorem 4.2]{bejamini-yadin-cylinder}}
For the external DLA process on $\G\times \mathbb{N}$, where $\G$ is a $d$-regular graph of size $n$, we have
$$D=\lim_{m\to\infty}\frac{1}{mn}\mathbb{E}[T_m]$$
\end{theorem}
In \cite[Theorem 4.6]{bejamini-yadin-cylinder} the previous result has been improved to $D\leq \frac{2}{3}$ for the case when the base graph $\G$ is a vertex transitive graph. Finally, for a family of base graphs with small mixing time, the following holds.
\begin{theorem}{\cite[Theorem 4.8]{bejamini-yadin-cylinder}}
Let $(\mathcal{G}^n)$ be a family of $d$-regular graphs ($d\geq 2$) such that $\lim_{n\to\infty}|\mathcal{G}^n|=\infty$, and for all $n$,
$$t_{\text{mix}}(\mathcal{G}^n)\leq \dfrac{\log^2|\mathcal{G}^n|}{(\log\log |\mathcal{G}^n|)^5}.$$
Let $D(n)$ be the density at infinity of the external DLA process on $\mathcal{G}^n\times \mathbb{N}$. Then
$\lim_{n\to\infty}D(n)=0$.
\end{theorem}
We refer to the last section of \cite{bejamini-yadin-cylinder} for several open questions and problems concerning external DLA on cylinder graphs. Many of the bounds from the previous three results can be improved, with some careful technicalities and assumptions on the base graph $\G$.

\subsection{Fractal graphs}\label{sec:frac-gr-dla}

\begin{small}
\begin{figure}
\label{fig:gasket}
	\centering
	\input{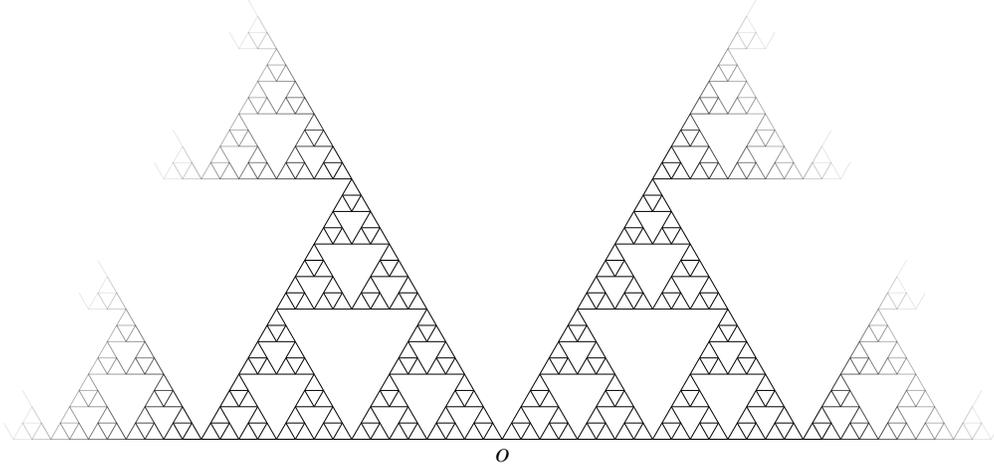}
	\caption{\label{fig:gasket} Doubly-infinite Sierpinski gasket graph $\SG$.}
\end{figure}
\end{small}

The appearance of fractal-like structures in DLA models (both internal and external) and their behavior on fractal graphs is the main theme of this paper, and we would like at this point to introduce two fractal graphs: the Sierpinski gasket graph and the Sierpinski carpet graph (called also pre-Sierpinki carpet).

\textbf{Sierpinski gasket graph $\SG$} is a pre-fractal associated with the Sierpinski gasket, defined as follows. We consider in $\mathbb{R}^2$ the sets
$V_0=\{(0,0), (1,0), (1/2,\sqrt{3}/2)\}$
and 
\begin{equation*}
E_0=\left\{\big((0,0),(1,0)\big),\big((0,0),(1/2,\sqrt{3}/2)\big),\big((1,0),(1/2,\sqrt{3}/2)\big)\right\}.
\end{equation*}
Now recursively define $(V_1,E_1), (V_2,E_2),\ldots$ by
\begin{equation*}
V_{n+1}=V_n\cup\left\{\big(2^n,0\big)+V_n\right\}\bigcup \left\{\left(2^{n-1},2^{n-1}\sqrt{3}\right)+V_n\right\}
\end{equation*}
and 
\begin{equation*}
E_{n+1}=E_n\cup\left\{\left(2^n,0\right)+E_n\right\}\bigcup \left\{\left(2^{n-1},2^{n-1}\sqrt{3}\right)+E_n\right\},
\end{equation*}
where $(x,y)+S:=\{(x,y)+s:s\in S\}$. Let $V_{\infty}=\cup_{n=0}^{\infty}V_n$, $E_{\infty}=\cup_{n=0}^{\infty}E_n$, $V=V_{\infty}\cup \{-V_{\infty}\}$ and $E=E_{\infty}\cup \{-E_{\infty}\}$. Then the doubly infinite Sierpinski gasket graph $\SG$
is the graph with vertex set $V$ and edge set $E$.
See Figure \ref{fig:gasket} for a graphical representation of $\SG$.
Set the origin $o=(0,0)$. External DLA on $\SG$ seems to be 
an approachable problem, due to the fact that $\SG$ is a post-critically finite fractal, and the existence of cut points simplifies the understanding of the harmonic measure from infinity, which can be defined again as in \eqref{eq:harm-measure-rec}, since the random walk on $\SG$ is recurrent.  We refer the reader to \cite{barlow-diffusion} and \cite{kigami-book} for more details on analysis and diffusion on fractals.

\textbf{Sierpinski carpet graph $\SC_m$,} called also $m$-dimensional pre-Sierpinski carpet, is an infinite graph derived from the \emph{Sierpinski carpet}. $\SC_2$ is constructed from  the unit square in $\R^2$ by dividing it 
into $9$ equal squares and deleting the one in the center. 
The same procedure is then repeated recursively to the remaining $8$ squares.
As mentioned in the introduction, we use the construction of the pre-Sierpinski carpet as in \cite{osada-pre-sierpinski}. Recall that in this construction, the length scale factor is $3$ and the mass scale factor is $3^m-1$.
 For random walks on such graphs see \cite{barlow_bass} and the references therein. See Figure \ref{fig:carpet} for a finite piece of Sierpinski carpet graph in
dimension $2$.

For $m\geq 3$, simple random walk on $\SC_m$ is transient, and the harmonic measure from infinity $\mu_A(y)$ for a finite subset $A\subset \SC_m$  is defined by using the capacity of $A$ and the equilibrium measure of $A$, similar to \eqref{eq:harm-measure-trans}. More details on the construction can be found in  \cite{benjamini-yadin-dla}, where upper bounds for the arms $r(\mathcal{E}_n)$ of external DLA on a large class of transient graphs, including $\SC_m$, $m\geq 3$, are proved. Their proofs are based on good control of heat-kernel estimates. The bounds for $\SC_m$ read as following.

\begin{figure}
\label{fig:carpet}
\begin{center}
\includegraphics[width=4.5cm]{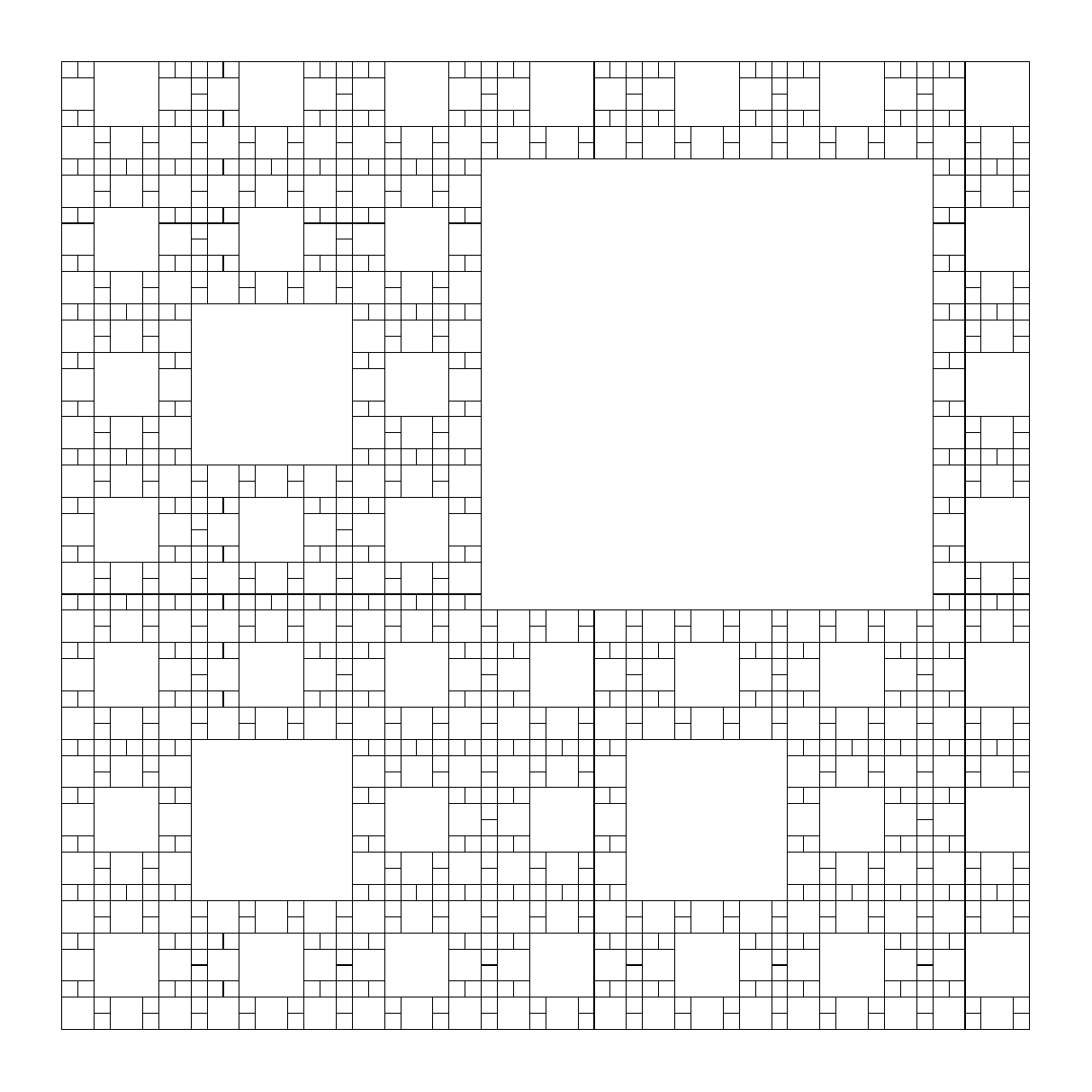}
\caption{\label{fig:carpet} Sierpinski carpet graph $\SC_2$}
\end{center}
\end{figure}
\begin{theorem}{\cite[Theorem 5.5]{benjamini-yadin-dla}}
Let $\SC_m$ be the $m$-dimensional Sierpinski carpet graph, and $(\mathcal{E}_n)_{n\geq 0}$ the external DLA process on $\SC_m$ started at $\mathcal{E}_0=\{o\}$ ($o$ is some fixed origin). Then almost surely,
$$\limsup_{n\to\infty} n^{-\beta} r(\mathcal{E}_n)<\infty$$
where
\begin{equation*}
\beta=
\begin{cases}
\frac{\log_2(13)-2}{3}=0.5568 &,\quad \text{if } m=3,\\
\frac{1}{2}&,\quad \text{if } m=4.
\end{cases}
\end{equation*}
When $m\geq 5$, we have almost surely,
$$\limsup_{n\to\infty}(\log n)^{-1}n^{-\frac{2}{d(m)-2}}r(\mathcal{E}_n)<\infty,$$
\end{theorem}
where $d(m)=\dfrac{\log(3^m-1)}{\log(3^m-1)-\log(3^{m-1}-1)}.$

We would like to conclude the section on external DLA with a couple of problems/questions.

\begin{question}
Can one find an upper bound for the growth of arms in external DLA on $\SG$ and on $\SC_2$ (the random walk is strongly recurrent on these two graphs)? Can one extend the method  Kesten used to upper bound the growth of arms in external DLA  on $\Z^2$?
\end{question}
\begin{question}
Do we have zero density at infinity of the cluster $\mathcal{E}_{\infty}$ on the Sierpinski gasket graph $\SG$?
\end{question}

\begin{question}
Does the external DLA cluster on the Sierpinski gasket graph and on the Sierpinski carpet graph have infinitely many holes, with probability one, as in the case of $\Z^2$ as proven in  \cite{eberz-wagner-dla}? 
\end{question}

Other than $\SG$ and $\SC_m$ there is a variety of other fractal graphs one can look at, and investigate the behavior of external DLA, which can be easier than $\Z^d$.

\begin{question}
Assuming that the Poisson boundary of the random walk on the graph $\G$ is non trivial, is there a characterization of the Poisson boundary in terms of the number of ends of the external DLA cluster at infinity $\mathcal{E}_{\infty}$ on $\G$?
\end{question}

\section{Internal DLA}\label{sec:idla}

Internal DLA can be defined on any infinite graph $\G$; fix as above a vertex $o$ of $\G$ and call it the origin.  The internal DLA cluster is built  up one site at a time, by letting the
$n$-th particle perform a random walk until it exits the set of sites already occupied by the previous $n-1$  particles, the walk of the $n$-th particle being independent of the past.
Similarly to external DLA, internal DLA is also a Markov chain on finite subsets of $\G$.

\begin{definition}\label{def:idla-mc} 
Let $\G$ be an infinite graph, and $(S_n)$ a simple random walk on $\G$ starting at  $o$. \emph{Internal DLA} on $\G$ is a Markov chain $(\mathcal{I}_n)_{n\geq 0}$ on finite connected subsets of $\G$, which evolves in the following way. Start with a single vertex $o\in \G$ and set $\mathcal{I}_0=\{o\}$. Given the process $\mathcal{I}_n$
at time $n$, let $y_{n+1}$ be a random vertex in $\partial\mathcal{I}_n$ chosen according to the harmonic measure (from $o$) of $\partial\mathcal{I}_n$, as defined in \eqref{eq:harm-meas}. That is, $y_{n+1}$ is the first exist location from $\mathcal{I}_n$ of the simple random walk $(S_n)$ starting from $o$, independent of the past:
$$\mathbb{P}[y_{n+1}=y|\mathcal{I}_n]=h_{\partial\mathcal{I}_n}(y),\quad \text{for }y\in\partial \mathcal{I}_n,$$ 
and we set $\mathcal{I}_{n+1}=\mathcal{I}_{n}\cup \{y_{n+1}\}.$ 
\end{definition}
The set $\mathcal{I}_{n}$ is called the \emph{internal DLA} cluster at time $n$, and it contains $n+1$ particles.
As $n\to\infty$, we are interested in the asymptotic shape of internal DLA cluster $\mathcal{I}_n$, and the fluctuations of the cluster around the limiting shape. Due to the fact that the harmonic measure for "nice subsets" (for example balls) of $\G$, when $\G$ is an Euclidean lattice, or a regular tree, is easier to understand than the harmonic measure from infinity as in the external DLA case, for the internal DLA model we have very precise estimates on many state spaces. Moreover, several variations of the classical internal DLA have been introduced.  

\subsection{Integer lattices $\mathbb{Z}^d$}

The first result concerning the internal DLA goes back to \cite{lawler_bramson_griffeath}, where it is shown that the limit shape of internal DLA cluster is a ball, in the following sense. Let $\omega_d$ be the volume of the $d$-dimensional Euclidean ball of radius 1, and $\mathcal{B}_n$ be the $d$-dimensional "lattice ball" of radius $n$, that is, $\mathcal{B}_n=\{x\in\Z^d:\ |x|\leq n\}$, where $|x|$ denotes the Euclidean norm of $x$. 
\begin{theorem}{\cite[Theorem 1]{lawler_bramson_griffeath}}
At time $\lfloor \omega_dn^d\rfloor$, internal DLA cluster occupies a set of sites close to a $d$-dimensional ball of radius $n$. More precisely, for any $\epsilon>0$, with probability $1$
$$\mathcal{B}_{n(1-\epsilon)}\subset \mathcal{I}_{\lfloor \omega_dn^d\rfloor}\subset \mathcal{B}_{n(1+\epsilon)}, \quad \text{ for n large}.$$
\end{theorem}
In this first paper, a basic open question on fluctuations (deviation of $\mathcal{I}_n$ from the Euclidean ball) was asked: are the fluctuations of order $\sqrt{n}$, of order $n^\delta$ for some $\delta\in(0,\frac{1}{2})$, or even smaller?  \textsc{Lawler} \cite{lawler_1995} proved that for $d\geq 2$, the fluctuations are subdiffusive and they are of order at most $n^{1/3}$. While it was conjectured that the fluctuations are at most logarithmic in the radius, this resisted a mathematical proof for about 20 years. Two independent groups \textsc{Jerison, Levine, and  Sheffield} \cite{jerison_levine_sheffield,jerison_levine_sheffield2,jerison_levine_sheffield3}
and \textsc{Asselah and Gaudilli\`{e}re} \cite{asselah_gaudilliere,asselah_gaudilliere_2,asselah_gaudilliere_3}, and by different methods have shown that indeed, for $d=2$ there are $\log(n)$ fluctuations, and for $d\geq 3$, there are $\sqrt{\log(n)}$ fluctuations in the radius. A summary of their results reads as following.
\begin{theorem}
If $d=2$, there is an absolute constant $c$, such that with probability $1$,
$$\mathcal{B}_{n- c\log n}\subset \mathcal{I}_{\lfloor \pi n^2\rfloor}\subset \mathcal{B}_{n+c\log n}, \text{ for all sufficiently large } n.$$
If $d\geq 3$, there is an absolute constant $C$, such that with probability $1$,
$$\mathcal{B}_{n- C\sqrt{\log n}}\subset \mathcal{I}_{\lfloor\omega_dn^d\rfloor}\subset \mathcal{B}_{n+C\sqrt{\log n}}, \text{ for all sufficiently large } n.$$
\end{theorem}
A generalization of the classical internal DLA on $\Z^d$ was treated in \cite{peres_levine_scaling_limits}, where instead of running all particles from the origin, 
the authors run the process from an arbitrary starting configuration of particles (initial density of particles) on finer and finer lattices, all particles still performing simple random walks. They then show that, as the lattice spacing tends to zero, the internal DLA has a deterministic scaling limit which can be described as the solution to a certain PDE free boundary problem in $\mathbb{R}^d$. We do not state here the rigorous result, which requires more notation and definition, but refer to the lengthy and complex paper \cite{peres_levine_scaling_limits}.
In order to study this general model, a new model called \emph{divisible sandpile} was introduced in \cite{peres_levine_strong_spherical}, which uses a continuous amount of mass instead of discrete particles. 

The \emph{divisible sandpile model} can be briefly described as following: start with an initial mass $\mu$ at the origin $o$. A vertex is called  \emph{full} if it has mass at least $1$. Any
full site can topple by keeping mass $1$ for itself and distributing the excess
mass equally among its neighbors. At each time step, one chooses a full site
and topples it. As time goes to infinity, provided each full site is eventually
toppled, the mass approaches a limiting distribution in which each site has mass $\leq 1$; this is proved in \cite{peres_levine_strong_spherical}. Individual topplings do not
commute, but the divisible sandpile is \emph{abelian} in the sense that any
sequence of topplings produces the same limiting mass distribution; this is
proved in \cite[Lemma 3.1]{peres_levine_scaling_limits}. The set of sites with limit mass distribution equal to $1$ is denoted by  $\mathcal{S}_n$ and is called \emph{the divisible sandpile cluster}.
The asymptotic shape of the divisible sandpile cluster $\mathcal{S}_n$ is proven to be the same as the one of the internal DLA cluster on $\Z^d$ in \cite{peres_levine_strong_spherical}, on regular trees in \cite{levine-sandpile-tree}, on comb lattices in \cite{huss_sava_aggregation}, and on Sierpinski gasket graphs in \cite{huss-sava-sandpile-gasket}.

\paragraph*{Random walks with drift on $\Z^d$.} If one lets the particles which build up the internal DLA cluster $\mathcal{I}_n$ perform drifted random walk instead of simple random walk as in the classical model, one can again ask  about the shape of the limit cluster on any state space. On $\Z^d$, this was open for several years, and the cluster was believed to be represented by the level sets of the Green function for the drifted random walk. This fact has been disproved, and with the help of the divisible sandpile model, in \cite{lucas-drift-idla} it was proven that the internal DLA cluster is a true heat ball, because it gives rise to a mean-value property for caloric functions. The author introduced there the \emph{unfair divisible sandpile}, where the mass is not distributed equally to the neighbors, but according to the one-step probabilities of the drifted random walk; the limit shape for the unfair divisible sandpile on $\Z^d$ was also described there. The main result for the limit shape for drifted internal DLA can be found in \cite[Theorem 1.1]{lucas-drift-idla}, and for the limit shape of the unfair divisible sandpile cluster in  \cite[Theorem 3.3]{lucas-drift-idla}. 

\paragraph*{Uniform starting points.} To my knowledge, the most recent result for internal DLA on $\Z^d$, concerns the limit shape for the cluster when the particles do not all start 
from the same vertex $o$. Instead the starting position is  chosen uniformly at random in the existing cluster. Formally, one can define the internal DLA as in Definition \ref{def:idla-mc}, starting with $\mathcal{I}_0=\{0\}$, and given the process at time $n$, let $y_{n+1}$ be the first exit location from $\mathcal{I}_n$ of the simple random walk $S_n^{X_n}$ starting at $X_n$, where $X_n$ is a point chosen uniformly on $\mathcal{I}_n$, independent of the past. Set $\mathcal{I}_{n+1}=\mathcal{I}_n\cup \{y_{n+1}\}$.
It turns out, as shown in \cite{uniform-dla}, that this additional source of randomness arising from the choice of the initial position of the random walk, does not change the limit shape of the process, as the result below shows. Let $b_n:=|\mathcal{B}_n|$.
\begin{theorem}{\cite[Theorem 1.1]{uniform-dla}}
Let $d\geq 2$. There exist constants $c_1,\ c_2,\ C_1$ and $C_2$ depending only on the dimension $d$ such that, almost surely, the internal DLA cluster $\mathcal{I}_n$ with uniform starting points satisfies 
$$\mathcal{B}_{n(1-C_1n^{-c_1})}\subseteq \mathcal{I}_{b_n}\subseteq \mathcal{B}_{n(1+C_2n^{-c_2})},\quad
\text{ for } n \text{ large enough}.$$
\end{theorem}

\begin{question}
What can we say about the fluctuations of the internal DLA cluster with uniform starting points around the limit shape? Are they bigger (smaller) that the fluctuations for internal DLA when all particles start their random walk from the same vertex $o$?
\end{question}

\textbf{Supercritical percolation cluster on $\Z^d$:} In \cite{shellef_idla_percolation}, the underlying state space for the internal DLA model is the supercritical bond percolation cluster on $\Z^d$, with the origin conditioned to be in the infinite cluster. It is shown in 
\cite[Theorem 1.1]{shellef_idla_percolation} that an inner bound for the internal DLA cluster is a ball in the graph metric. The picture for the outer bound was completed in \cite[Theorem 1.1]{duminil_lucas_yadin_yehudayoff}, where the authors show that also in this case the limit shape is a ball. The results in their paper hold in a more general setting: given the existence of a "good" inner bound for internal DLA, one can also prove a matching outer bound by using their methods. An interesting problem in the context of internal DLA model on a random graph is to understand the fluctuations.

\subsection{Comb lattices $\mathcal{C}_2$}

The 2-dimensional comb lattice $\C$ is the spanning tree of $\Z^2$ obtained by removing all horizontal edges, except the ones on the $x$-axis. While $\C$ is a simple graph, see Figure \ref{fig:comb}, it has some remarkable properties in what concerns the behavior of random walks: no form of the so-called Einstein relation for exponents associated
with random walks hold on $\C$, see \cite{bertacchi_comb}. \textsc{Peres and Krishnapur} \cite{peres_krishnapur_collide} showed that on $\C$ two independent simple random walks meet only finitely often. The comb $\C$ is an example where the limit shape of internal DLA is not a ball in the graph metric or in another standard metric.
Indeed, the diameter of the internal DLA cluster with $n$ particles grows like $n^{2/3}$ in the $y$-direction, and like $n^{1/3}$ in the $x$-direction. See Figure \ref{fig:idla_comb} for a picture of the internal DLA cluster with 100, 500, and 1000 particles, respectively.

\begin{minipage}{0.40\linewidth}
\begin{center}
\begin{tikzpicture}[scale=1.1]
\foreach \x in {-2,...,2}
{
   \draw (\x, 2.5) -- (\x, -2.5);
   \foreach \y in {-2,...,2}
   \fill (\x,\y) circle (1.5pt);
}

\draw (-2.5, 0) -- (2.5, 0);
\coordinate[label=135:$o$] (O) at (0,0);
\end{tikzpicture}
\captionof{figure}{\label{fig:comb}The comb $\C$.}
\end{center}
\end{minipage}
\hfill
\begin{minipage}{0.60\linewidth}
\begin{center}
  \begin{tikzpicture}[scale=0.06]

\begin{scope}[xshift=-20cm]
\draw (-4,0) -- (-4,1);
\draw (-3,-3) -- (-3,3);
\draw (-2,-5) -- (-2,4);
\draw (-1,-9) -- (-1,8);
\draw (0,-11) -- (0,12);
\draw (1,-9) -- (1,10);
\draw (2,-4) -- (2,5);
\draw (3,-1) -- (3,2);
\draw (4,-1) -- (4,1);
\draw (5,0) -- (5,1);
\draw (6,0) -- (6,0);
\draw (-4,0) -- (6,0);
\end{scope}

\begin{scope}
\draw (-10,0) -- (-10,0);
\draw (-9,-2) -- (-9,1);
\draw (-8,-3) -- (-8,3);
\draw (-7,-2) -- (-7,5);
\draw (-6,-5) -- (-6,8);
\draw (-5,-11) -- (-5,8);
\draw (-4,-14) -- (-4,14);
\draw (-3,-17) -- (-3,19);
\draw (-2,-24) -- (-2,24);
\draw (-1,-26) -- (-1,26);
\draw (0,-40) -- (0,33);
\draw (1,-26) -- (1,29);
\draw (2,-19) -- (2,21);
\draw (3,-18) -- (3,18);
\draw (4,-12) -- (4,16);
\draw (5,-10) -- (5,8);
\draw (6,-4) -- (6,5);
\draw (7,-1) -- (7,2);
\draw (8,-1) -- (8,3);
\draw (9,-1) -- (9,2);
\draw (-10,0) -- (9,0);
\end{scope}

\begin{scope}[xshift=26cm]
\draw (-12,0) -- (-12,0);
\draw (-11,-1) -- (-11,1);
\draw (-10,-2) -- (-10,5);
\draw (-9,-5) -- (-9,5);
\draw (-8,-8) -- (-8,9);
\draw (-7,-10) -- (-7,12);
\draw (-6,-13) -- (-6,15);
\draw (-5,-23) -- (-5,24);
\draw (-4,-24) -- (-4,30);
\draw (-3,-31) -- (-3,33);
\draw (-2,-42) -- (-2,41);
\draw (-1,-45) -- (-1,47);
\draw (0,-60) -- (0,51);
\draw (1,-52) -- (1,48);
\draw (2,-35) -- (2,45);
\draw (3,-38) -- (3,33);
\draw (4,-25) -- (4,28);
\draw (5,-27) -- (5,22);
\draw (6,-13) -- (6,15);
\draw (7,-9) -- (7,12);
\draw (8,-9) -- (8,7);
\draw (9,-7) -- (9,5);
\draw (10,-3) -- (10,4);
\draw (11,0) -- (11,2);
\draw (12,0) -- (12,0);
\draw (-12,0) -- (12,0);
\end{scope}

\end{tikzpicture}
    
  
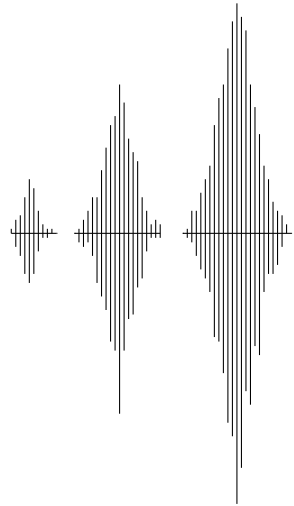
\captionof{figure}{\label{fig:idla_comb}Internal DLA cluster.}
  \end{center}
\end{minipage}

Let 
\begin{align}\label{eq:limit_shape}
\mathcal{D}_n = \left\lbrace (x,y)\in\mathcal{C}_2:\: \frac{\abs{x}}{k} + 
\left(\frac{\abs{y}}{l}\right)^{1/2}\leq n^{1/3}\right\rbrace
\end{align}
where the constants $k$ and $l$ are given by
\begin{equation*}
k = \left(\frac{3}{2}\right)^{2/3},
\quad l = \frac{1}{2}\left(\frac{3}{2}\right)^{1/3}.
\end{equation*}
The  inner bound for the limit shape of internal DLA cluster on $\C$ was proven in \cite[Theorem 4.2]{huss_sava_aggregation}, while the outer bound together together with the fluctuations was proven in \cite{asselah_comb}. 

\begin{theorem}{\cite[Theorem 1.2]{asselah_comb}}
There is a positive constant $a$ such that with probability $1$, and $n$ large enough
$$\mathcal{D}_{n-a\sqrt{\log n}}\subset \mathcal{I}_n\subset \mathcal{D}_{n+a\sqrt{\log n}}.$$
\end{theorem}
Remark that this result does not mean that the fluctuations are sub-logarithmic, but rather gaussian; see \cite[Theorem 1.2]{asselah_comb} and the comments afterwards. In \cite[Theorem 3.5]{huss_sava_aggregation} we also prove that the limit shape for the divisible sandpile cluster on $\C$ is given by the set $\mathcal{D}_n$.

\subsection{Trees $\mathbb{T}_d$}

Internal DLA on discrete groups with exponential growth has been studied in  \cite{blachere_brofferio}. The homogeneous tree $\mathbb{T}_d$ is a particular case (as a Cayley graph of a free group) of these state spaces, for which the authors have proven that the limit shape of internal DLA cluster is a ball in the graph metric, and they give  lower bounds for the inner and outer error. The more general result is the following.
\begin{theorem}\cite[Theorem 3.1]{blachere_brofferio}
Let $\G$ be a finitely generated group of exponential growth, and consider the internal DLA model $(\mathcal{I}_n)$ on $\G$, built up with symmetric random walks with finitely supported increments, starting at the identity $o$ of $\G$. Then, for any constants $C_O>2$ and $C_I>3/K$,
$$\mathbb{P}\big[\exists n_0\ s.t.\ \forall n>n_0:\quad \B_{n-C_I\ln n}\subset \mathcal{I}_{|\B_n|}\subset \B_{n+C_O\sqrt{n}}\big]=1,$$
where $K$ is a constant that ensures that the ball $\B_n$ contains the boundary $\partial \B_{n-1}$,
and $\B_n$ is the ball of radius $n$ centered at the identity in the word metric on $\G$.
\end{theorem}
An extension of this result to non-amenable graphs for a wide class of Markov chains was considered in \cite{huss_idla}. On discrete groups with polynomial growth, internal DLA has been considered in \cite{blachere_2004}.

\subsection{Cylinder graphs}

Like in Section \ref{sec:cyl-gr-dla}, we consider here cylinder graphs $\G\times \mathbb{Z}$, and we let $\G$ to be the cycle graph $\mathbb{Z}_N$ on $N$ vertices. Internal DLA on cylinder graphs $\Z_N\times \mathbb{Z}$ was investigated in \cite{jerison_levine_sheffield-cylinders}, for the following initial setting. For $k\in \Z$, the set $\Z_{N}\times \{k\}$  is called the $k$-th level of the cylinder, and $R_k=\{(x,y)\in \Z_N\times \mathbb{N}:\ y\leq k\}$ the \emph{rectangle of height } $k$.
Let $\mathcal{I}_0=R_0$, and given the cluster $\mathcal{I}_n$ at time $n$, let $y_{n+1}$ be the first exit location from $\mathcal{I}_n$ of a random walk that starts uniformly at random on level zero of the cylinder, independent on the past, that is, the starting location is chosen with equal probability among the $N$ sites $(x,0)$, for $x\in\mathbb{Z}_N$. We then set $\mathcal{I}_{n+1}=\mathcal{I}_n\cup \{y_{n+1}\}$. It has been proven in \cite[Theorem 2]{jerison_levine_sheffield-cylinders} that the limit shape of internal DLA clusters on $\Z_N\times \mathbb{Z}$ is logarithmically close to rectangles, result that we do not state in complete form here, but instead we state a more recent result due to \textsc{Levine and Silvestri} \cite[Theorem 1.1]{levine_silvestri-cylinder} which generalizes the previous one \cite{jerison_levine_sheffield-cylinders} (here  the fluctuations are described in terms of the Gaussian Free Field exactly). Remark that in the cylinder setting, there are two parameters, the size $N$ of the cycle base graph, and the time $n$.
\begin{theorem}{\cite[Theorem 1.1]{levine_silvestri-cylinder}}
Let $(\mathcal{I}_n)_{n\geq 0}$ be the internal DLA process on $\Z_N\times \mathbb{Z}$ starting from $\mathcal{I}_0=R_0$. For any $\gamma >0$, $m\in\mathbb{N}$ there exist a constant $C=C(\gamma,m)$ such that
$$\mathbb{P}[R_{\frac{n}{N}-C\log N}\subseteq \mathcal{I}_n\subseteq R_{\frac{n}{N}+C\log N},\ n\leq N^m ]\geq 1-N^{-\gamma},\ \text{for N large enough}.$$

\end{theorem}
For other results concerning the fluctuations and the behavior of internal DLA clusters on $\Z_N\times \mathbb{Z}$, we refer to \cite{levine_silvestri-cylinder}.

\subsection{Fractal graphs}

We would like to conclude the section about internal DLA with the behavior of the model on Sierpinski gasket graphs $\SG$. Recall the definition of the Sierpinski gasket graph $\SG$ and of the Sierpinski carpet graph $\SC_2$, as given in Section \ref{sec:frac-gr-dla}. Due to the symmetry of $\SG$, it is clear that the limit shape of the internal DLA cluster on $\SG$ is a ball in the  graph metric, a result proved in \cite{chen-huss-sava-teplyaev}.

\begin{figure}[H]
\centering
\begin{tikzpicture}
\node (a1)          at (0,0) {\scalebox{1}[-1]{\includegraphics[width=3.2cm]{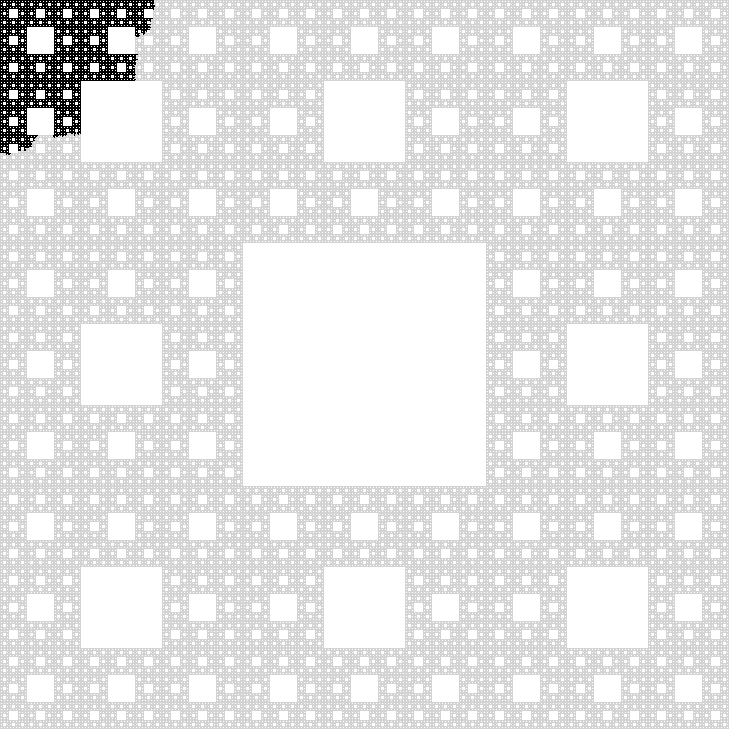}}};
\node (b1) [right=0.5cm of a1] {\scalebox{1}[-1]{\includegraphics[width=3.2cm]{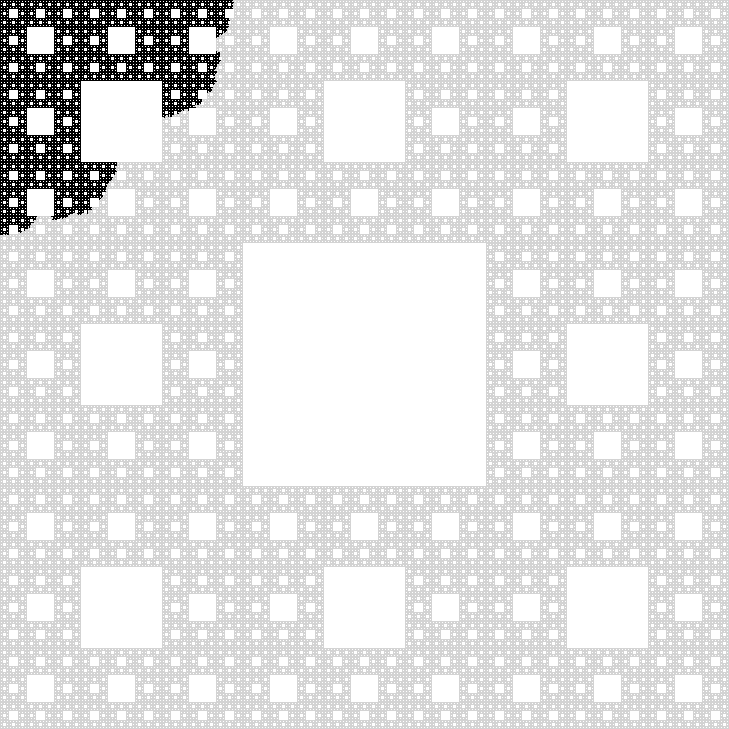}}};
\node (c1) [right=0.5cm of b1] {\scalebox{1}[-1]{\includegraphics[width=3.2cm]{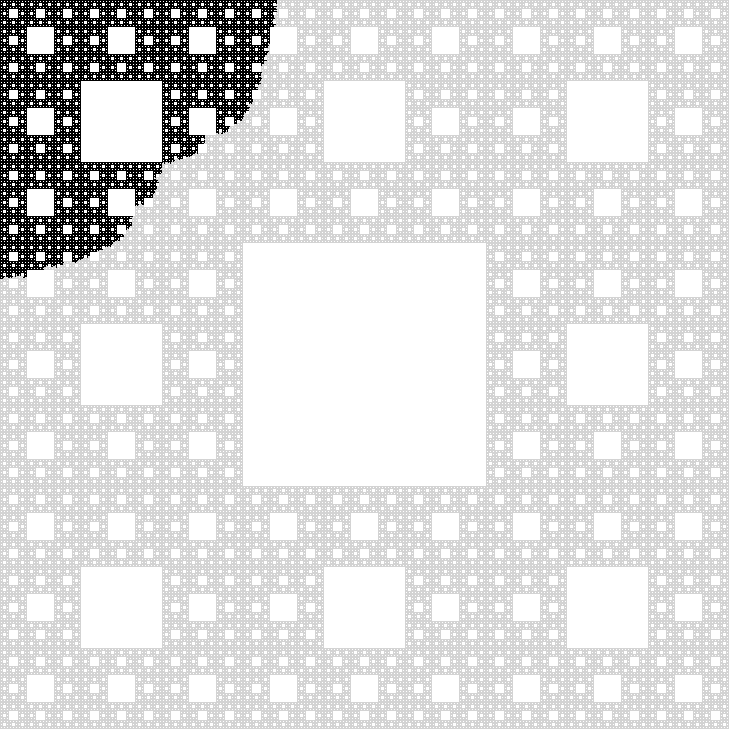}}};
\node (d1) [right=0.5cm of c1] {\scalebox{1}[-1]{\includegraphics[width=3.2cm]{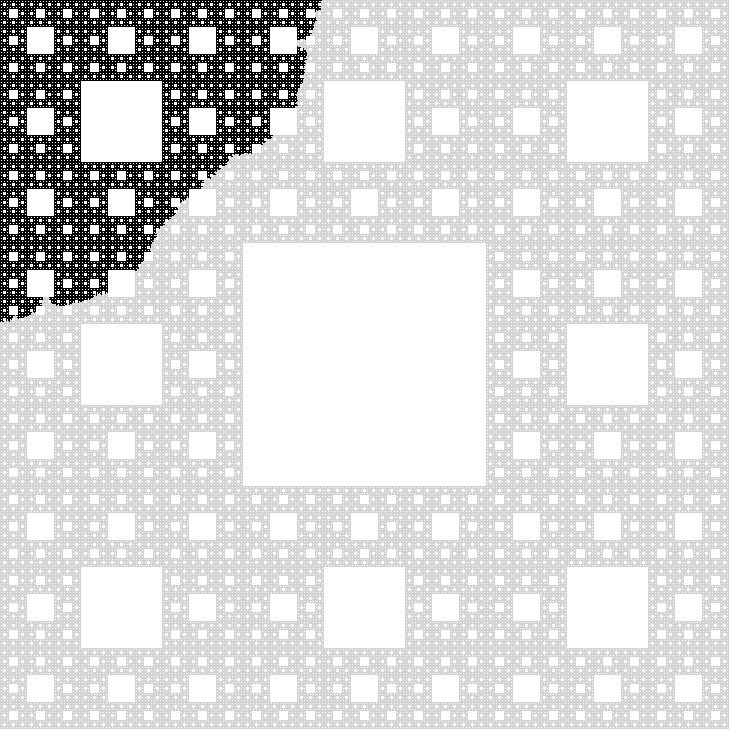}}};
\node (a2) [below=0.5cm of a1] {\scalebox{1}[-1]{\includegraphics[width=3.2cm]{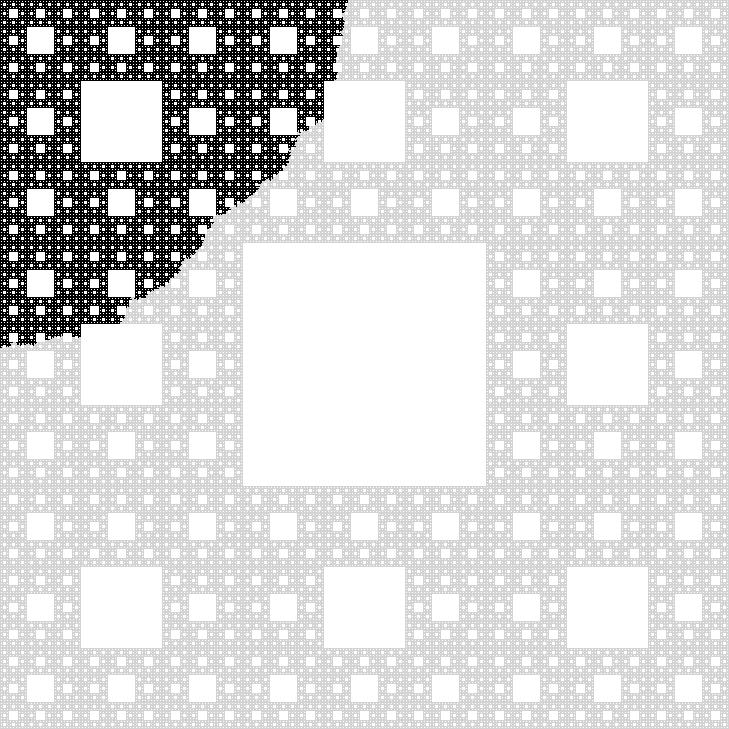}}};
\node (b2) [right=0.5cm of a2] {\scalebox{1}[-1]{\includegraphics[width=3.2cm]{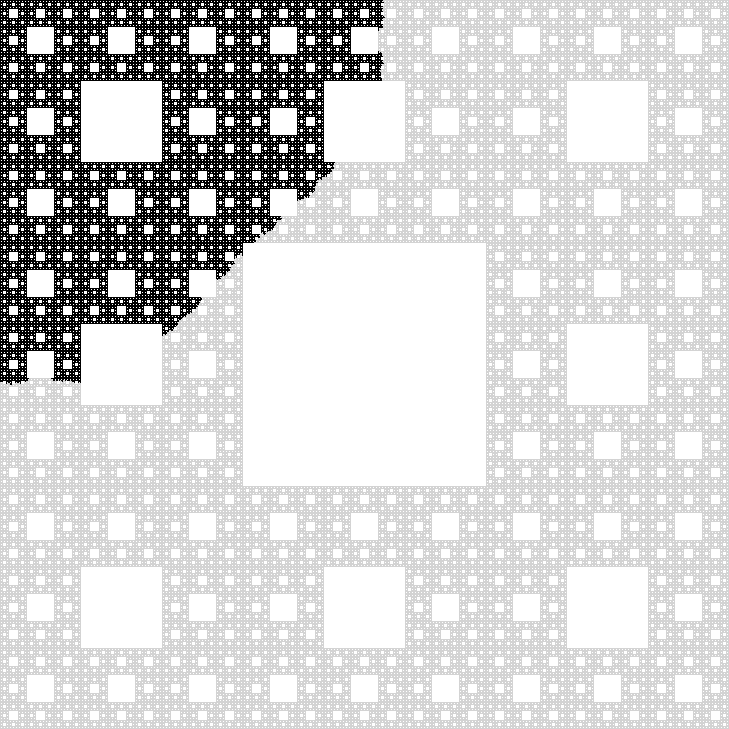}}};
\node (c2) [right=0.5cm of b2] {\scalebox{1}[-1]{\includegraphics[width=3.2cm]{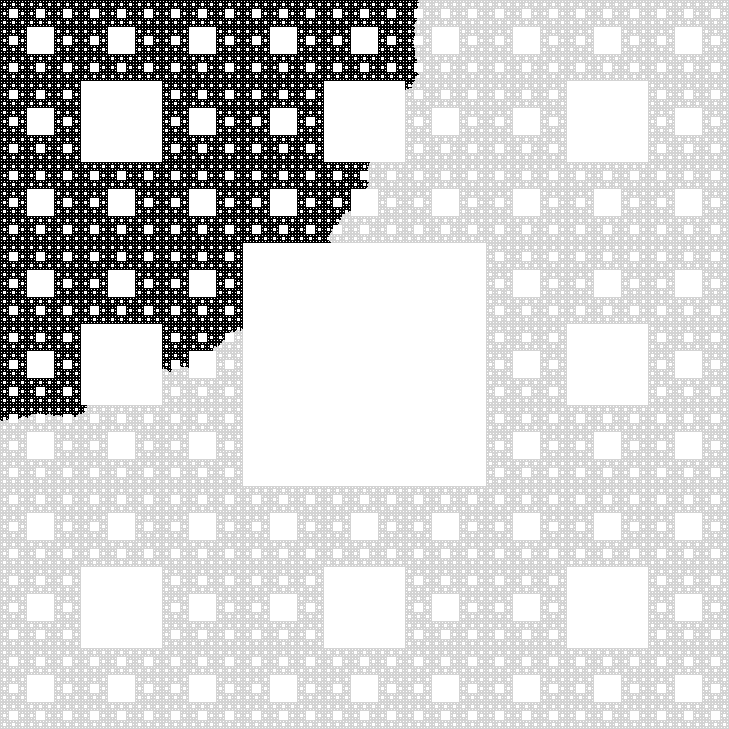}}};
\node (d2) [right=0.5cm of c2] {\scalebox{1}[-1]{\includegraphics[width=3.2cm]{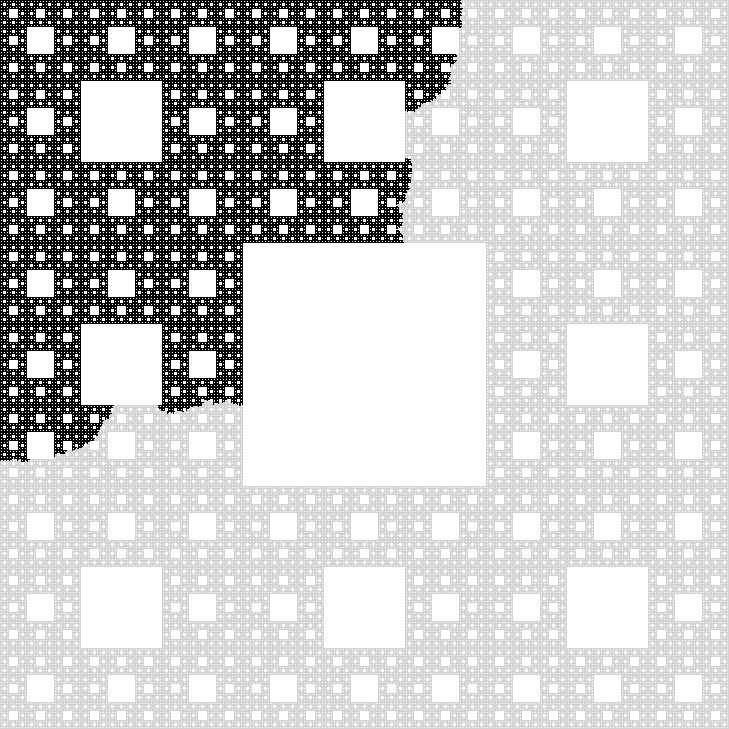}}};
\node (a3) [below=0.5cm of a2] {\scalebox{1}[-1]{\includegraphics[width=3.2cm]{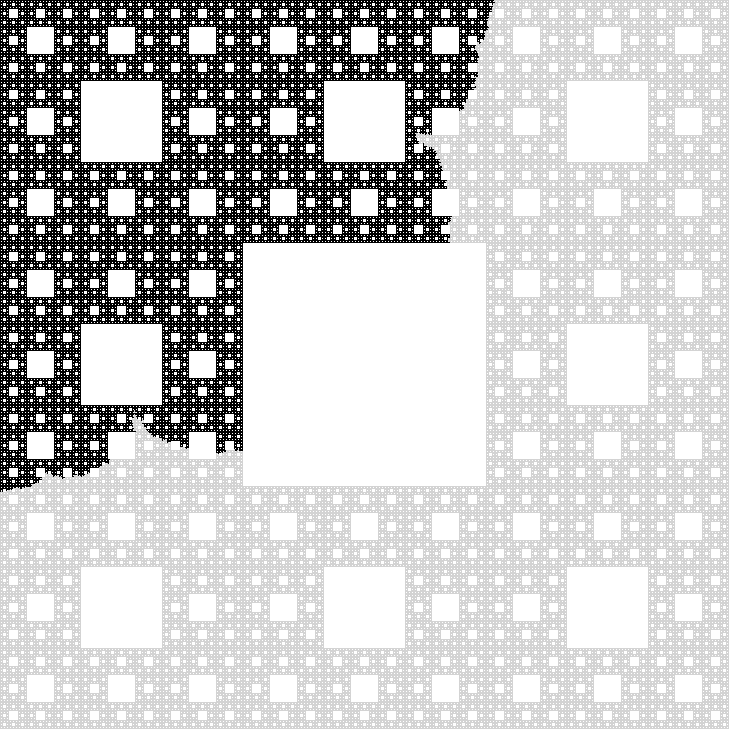}}};
\node (b3) [right=0.5cm of a3] {\scalebox{1}[-1]{\includegraphics[width=3.2cm]{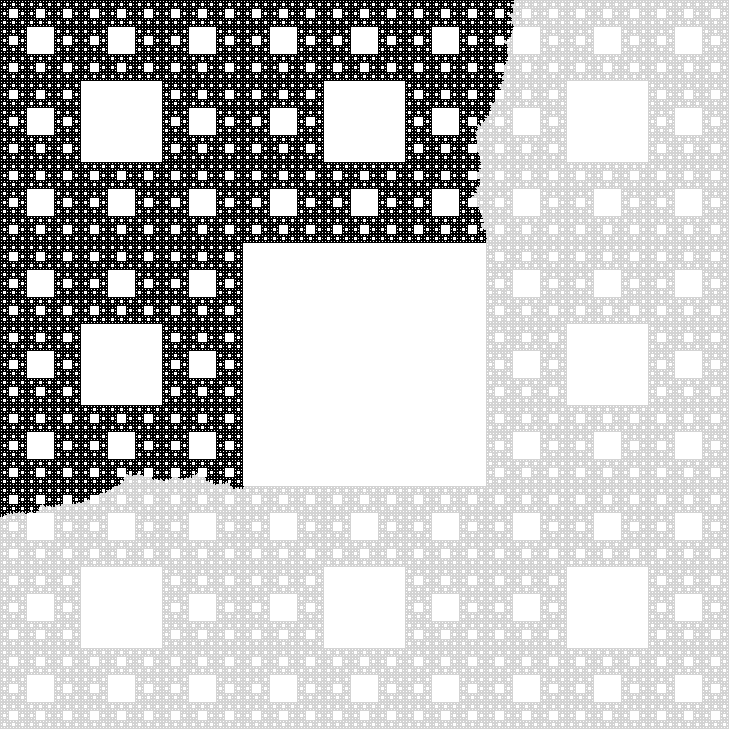}}};
\node (c3) [right=0.5cm of b3] {\scalebox{1}[-1]{\includegraphics[width=3.2cm]{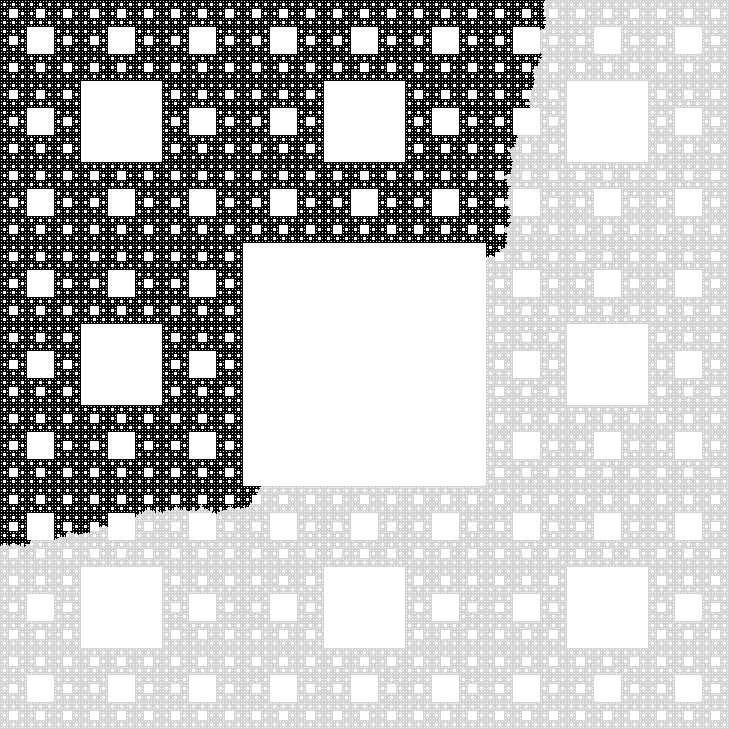}}};
\node (d3) [right=0.5cm of c3] {\scalebox{1}[-1]{\includegraphics[width=3.2cm]{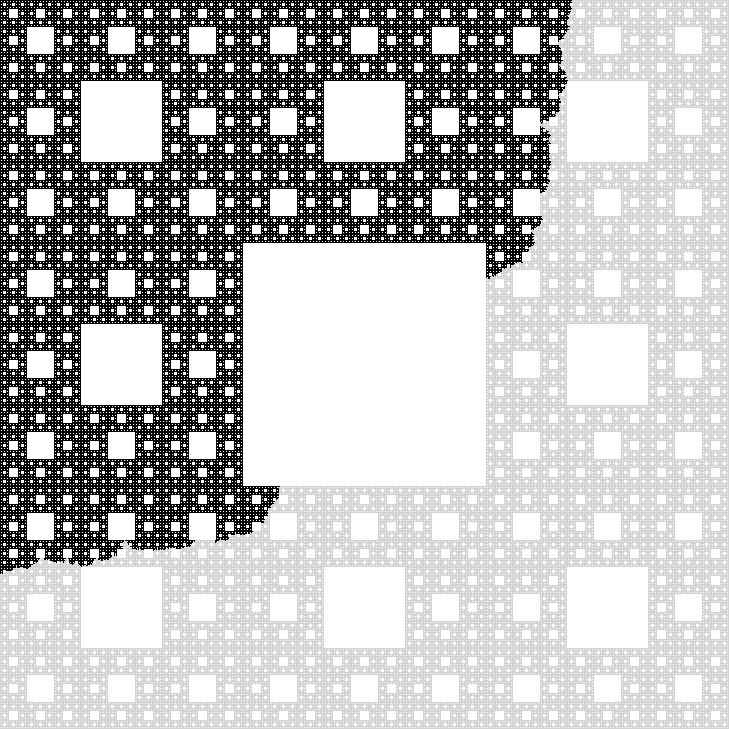}}};
\node (a4) [below=0.5cm of a3] {\scalebox{1}[-1]{\includegraphics[width=3.2cm]{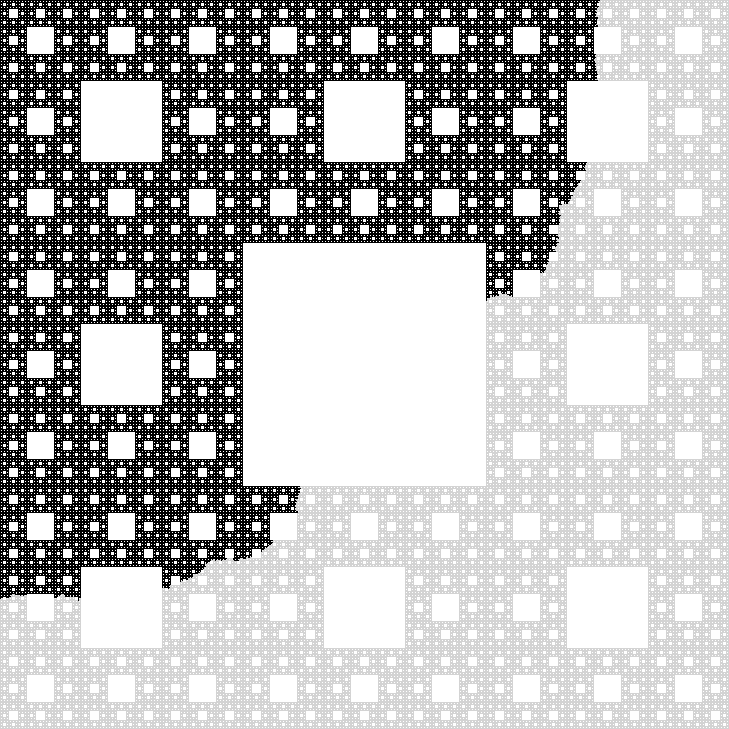}}};
\node (b4) [right=0.5cm of a4] {\scalebox{1}[-1]{\includegraphics[width=3.2cm]{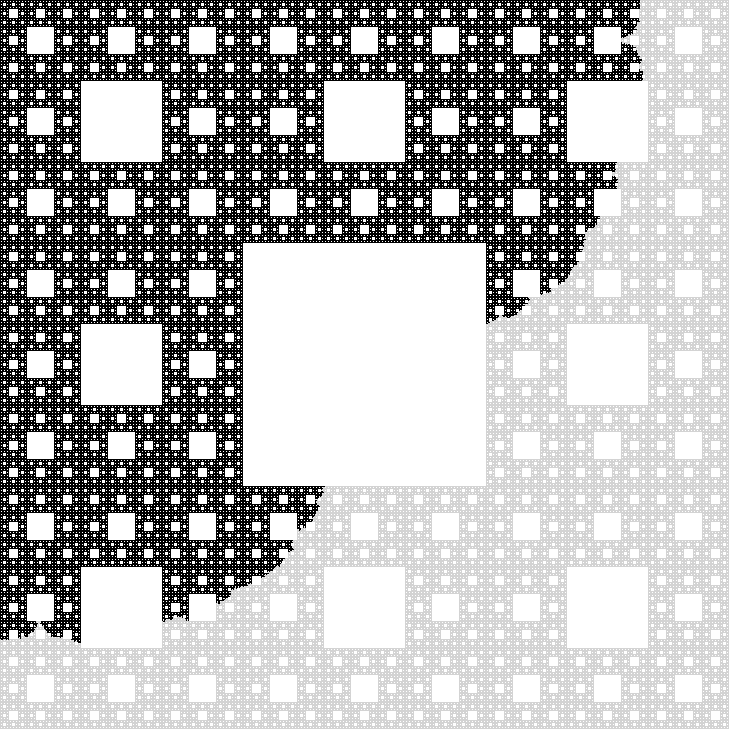}}};
\node (c4) [right=0.5cm of b4] {\scalebox{1}[-1]{\includegraphics[width=3.2cm]{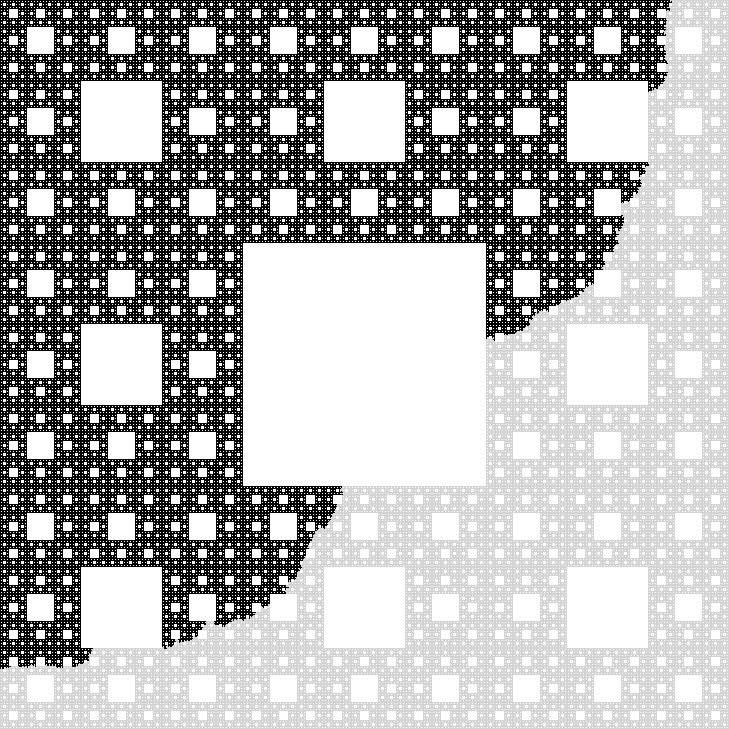}}};
\end{tikzpicture}
\caption{\label{fig:idla_sc}Internal DLA clusters on $\SC_2$ for $10000$ up to $150000$ particles.\\ Simulations by W.Huss}
\end{figure}

\begin{theorem}{\cite[Theorem 1.1]{chen-huss-sava-teplyaev}}
On $\SG$, the internal DLA cluster of $|\mathcal{B}_n|$ particles occupies a set of sites close to
a ball of radius $n$. That is, for all $\epsilon >0 $, we have
\begin{equation*}
\mathcal{B}_{n(1-\epsilon)}\subset \mathcal{I}_{|\mathcal{B}_n|}\subset \mathcal{B}_{n(1+\epsilon)}, \text{ for all n sufficiently large}
\end{equation*}
with probability 1.
\end{theorem}
A limit shape for the divisible sandpile on $\SG$ was described in \cite{huss-sava-sandpile-gasket}. Concerning the fluctuations for internal DLA, it is conjectured that they are sub-logarithmic.
\begin{conjecture}{\cite[Conjecture 4.1]{chen-gasket-universality}}
There exists $C>0$ such that
$$\mathcal{B}_{n-C\sqrt{\log n}}\subset \mathcal{I}_{|\mathcal{B}_n|}\subset \mathcal{B}_{n+C\sqrt{\log n}}.$$
\end{conjecture}
Many other questions concerning internal DLA on fractal graphs can be found in the final section of \cite{chen-gasket-universality}.
\begin{question}
Is the limit shape for the internal DLA model with uniform starting points on $\SG$, again a ball in the graph metric? What about the fluctuations in this case?
\end{question}
A reason why $\SG$ is easier to work with is because 1) it is a finitely ramified fractal graph, and 2) we have a precise characterization of the divisible sandpile model on $\SG$, thanks to the finite ramification and the symmetries it possesses. In contrast, $\SC_2$ is infinitely ramified, and characterizing the harmonic measure thereon is a challenging open question in the study of analysis on fractals. So at the moment it is very difficult to analyze growth models on $\SC_2$. See Figure \ref{fig:idla_sc} for the behavior of internal DLA on $\SC_2$.

\begin{question}
Does the internal DLA cluster on the 2-dimensional Sierpinski carpet graph $\SC_2$ have a (unique) scaling limit? What can one say about the boundary of the limit shape, which according to simulations appears to be of fractal nature?
\end{question}
\begin{question}
What is the limit shape of internal DLA on fractal graphs, other that $\SG$ (which is understood)
and $\SC_2$ (which seems hard to investigate)?
\end{question}
Since in most cases, the limit shape for internal DLA is a ball (in the graph metric, or Euclidean metric, or word metric), a more general question to ask is about the state space for the process.
\begin{question}
What properties should the state space $\G$ and the random walk on it have, in order for the internal DLA cluster on $\G$ to have a ball as limit shape? 
\end{question}

We would like to conclude this survey with the remark that fractals provide a class of state spaces with intriguing properties, both for the behavior of the external and internal DLA model, respectively. This behavior is definitely not fully understood on such graphs, and we hope to attract more people from the fractal community into the beauty of these topics.

\paragraph*{Acknowledgements.}
I am very grateful to the anonymous referee for a very careful reading of the manuscript and for several useful comments that improved the paper substantially.

\bibliography{proposal}{}

\begin{thebibliography}{DCLYY13}

\bibitem[AABK16]{dla-long-range1}
Gideon Amir, Omer Angel, Itai Benjamini, and Gady Kozma.
\newblock One-dimensional long-range diffusion-limited aggregation {I}.
\newblock {\em Ann. Probab.}, 44(5):3546--3579, 2016.

\bibitem[AAK17]{dla-long-range2}
Gideon Amir, Omer Angel, and Gady Kozma.
\newblock One-dimensional long-range diffusion limited aggregation {II}: {T}he
  transient case.
\newblock {\em Ann. Appl. Probab.}, 27(3):1886--1922, 2017.

\bibitem[AG13a]{asselah_gaudilliere_2}
Amine Asselah and Alexandre Gaudilli{\`e}re.
\newblock From logarithmic to subdiffusive polynomial fluctuations for internal
  {DLA} and related growth models.
\newblock {\em Ann. Probab.}, 41(3A):1115--1159, 2013.

\bibitem[AG13b]{asselah_gaudilliere_3}
Amine Asselah and Alexandre Gaudilli{\`e}re.
\newblock Sublogarithmic fluctuations for internal {DLA}.
\newblock {\em Ann. Probab.}, 41(3A):1160--1179, 2013.

\bibitem[AG14]{asselah_gaudilliere}
Amine Asselah and Alexandre Gaudilli{\`e}re.
\newblock Lower bounds on fluctuations for internal {DLA}.
\newblock {\em Probab. Theory Related Fields}, 158(1-2):39--53, 2014.

\bibitem[Ami17]{amir-dla}
Gideon Amir.
\newblock One-dimensional long-range diffusion-limited aggregation
  {III}---{T}he limit aggregate.
\newblock {\em Ann. Inst. Henri Poincar\'{e} Probab. Stat.}, 53(4):1513--1527,
  2017.

\bibitem[AR16]{asselah_comb}
Amine Asselah and Houda Rahmani.
\newblock Fluctuations for internal {DLA} on the comb.
\newblock {\em Ann. Inst. Henri Poincar\'{e} Probab. Stat.}, 52(1):58--83,
  2016.

\bibitem[Bar98]{barlow-diffusion}
Martin~T. Barlow.
\newblock Diffusions on fractals.
\newblock In {\em Lectures on probability theory and statistics
  ({S}aint-{F}lour, 1995)}, volume 1690 of {\em Lecture Notes in Math.}, pages
  1--121. Springer, Berlin, 1998.

\bibitem[BB99]{barlow_bass}
Martin~T. Barlow and Richard~F. Bass.
\newblock Random walks on graphical {S}ierpinski carpets.
\newblock In {\em Random walks and discrete potential theory ({C}ortona,
  1997)}, Sympos. Math., XXXIX, pages 26--55. Cambridge Univ. Press, Cambridge,
  1999.

\bibitem[BB07]{blachere_brofferio}
S{\'e}bastien Blach{\`e}re and Sara Brofferio.
\newblock Internal diffusion limited aggregation on discrete groups having
  exponential growth.
\newblock {\em Probab. Theory Related Fields}, 137(3-4):323--343, 2007.

\bibitem[BDCKL19]{uniform-dla}
Itai Benjamini, Hugo Duminil-Copin, Gady Kozma, and Cyrille Lucas.
\newblock Internal diffusion-limited aggregation with uniform starting points.
\newblock {\em Annales de l'Institut Henri Poincaré}, 2019,
  \href{http://arxiv.org/abs/1707.03241}{{\ttfamily arXiv:1707.03241}}.
\newblock to appear.

\bibitem[Ber06]{bertacchi_comb}
Daniela Bertacchi.
\newblock Asymptotic behaviour of the simple random walk on the 2-dimensional
  comb.
\newblock {\em Electron. J. Probab.}, 11:no. 45, 1184--1203, 2006.

\bibitem[Bla04]{blachere_2004}
S{\'e}bastien Blach{\`e}re.
\newblock Internal diffusion limited aggregation on discrete groups of
  polynomial growth.
\newblock In {\em Random walks and geometry}, pages 377--391. Walter de Gruyter
  GmbH \& Co. KG, Berlin, 2004.

\bibitem[BPP97]{barlow-pemantle-perkins-dla}
Martin~T. Barlow, Robin Pemantle, and Edwin~A. Perkins.
\newblock Diffusion-limited aggregation on a tree.
\newblock {\em Probab. Theory Related Fields}, 107(1):1--60, 1997.

\bibitem[BY08]{bejamini-yadin-cylinder}
Itai Benjamini and Ariel Yadin.
\newblock Diffusion limited aggregation on a cylinder.
\newblock {\em Comm. Math. Phys.}, 279(1):187--223, 2008.

\bibitem[BY17]{benjamini-yadin-dla}
Itai Benjamini and Ariel Yadin.
\newblock Upper bounds on the growth rate of diffusion limited aggregation.
\newblock 2017,  \href{http://arxiv.org/abs/1705.06095}{{\ttfamily
  arXiv:1705.06095}}.
\newblock preprint.

\bibitem[CHSHT19]{chen-huss-sava-teplyaev}
Joe~P. Chen, Wilfried Huss, Ecaterina Sava-Huss, and Alexander Teplyaev.
\newblock Internal {DLA} on {S}ierpinski gasket graphs.
\newblock {\em Analysis and Geometry on Graphs and Manifolds, London
  Mathematical Society Lecture Note Series}, 2019,
  \href{http://arxiv.org/abs/1702.04017}{{\ttfamily arXiv:1702.04017}}.
\newblock to appear.

\bibitem[CKF19]{chen-gasket-universality}
Joe~P. Chen and Jonah Kudler-Flam.
\newblock Laplacian growth and sandpiles on the {S}ierpinski gasket: limit
  shape universality and exact solutions.
\newblock {\em Ann. Inst.Henri Poincaré Comb. Phys. Interact.}, 2019.
\newblock to appear.

\bibitem[DCLYY13]{duminil_lucas_yadin_yehudayoff}
Hugo Duminil-Copin, Cyrille Lucas, Ariel Yadin, and Amir Yehudayoff.
\newblock Containing internal diffusion limited aggregation.
\newblock {\em Electron. Commun. Probab.}, 18:no. 50, 8, 2013.

\bibitem[DF91]{diaconis_fulton_1991}
Persi Diaconis and W.~Fulton.
\newblock A growth model, a game, an algebra, {L}agrange inversion, and
  characteristic classes.
\newblock {\em Rend. Sem. Mat. Univ. Politec. Torino}, 49(1):95--119 (1993),
  1991.
\newblock Commutative algebra and algebraic geometry, II (Italian) (Turin,
  1990).

\bibitem[Eld15]{eldan-dla-hyper}
Ronen Eldan.
\newblock Diffusion-limited aggregation on the hyperbolic plane.
\newblock {\em Ann. Probab.}, 43(4):2084--2118, 2015.

\bibitem[EW99]{eberz-wagner-dla}
Dorothea~M. Eberz-Wagner.
\newblock {\em Discrete Growth Models}.
\newblock PhD thesis, University of Washington, 1999.

\bibitem[HS12]{huss_sava_aggregation}
Wilfried Huss and Ecaterina Sava.
\newblock Internal aggregation models on comb lattices.
\newblock {\em Electron. J. Probab.}, 17:no. 30, 21, 2012.

\bibitem[HSH19]{huss-sava-sandpile-gasket}
Wilfried Huss and Ecaterina Sava-Huss.
\newblock Divisible sandpile on {S}ierpinski gasket graphs.
\newblock {\em Fractals}, 27, 2019,
  \href{http://arxiv.org/abs/1702.08370}{{\ttfamily arXiv:1702.08370}}.

\bibitem[Hus08]{huss_idla}
Wilfried Huss.
\newblock Internal {D}iffusion-{L}imited {A}ggregation on non-amenable graphs.
\newblock {\em Electronic Communications in Probability}, 13:272--279, 2008.

\bibitem[JLS12]{jerison_levine_sheffield}
David Jerison, Lionel Levine, and Scott Sheffield.
\newblock Logarithmic fluctuations for internal {DLA}.
\newblock {\em J. Amer. Math. Soc.}, 25(1):271--301, 2012.

\bibitem[JLS13]{jerison_levine_sheffield2}
David Jerison, Lionel Levine, and Scott Sheffield.
\newblock Internal {DLA} in higher dimensions.
\newblock {\em Electronic Journal of Probability}, 18(98):1--14, 2013,
  \href{http://arxiv.org/abs/1012.3453}{{\ttfamily arXiv:1012.3453}}.

\bibitem[JLS14a]{jerison_levine_sheffield3}
David Jerison, Lionel Levine, and Scott Sheffield.
\newblock Internal {DLA} and the {G}aussian free field.
\newblock {\em Duke Math. J.}, 163(2):267--308, 2014.

\bibitem[JLS14b]{jerison_levine_sheffield-cylinders}
David Jerison, Lionel Levine, and Scott Sheffield.
\newblock Internal {DLA} for cylinders.
\newblock In {\em Advances in analysis: the legacy of {E}lias {M}. {S}tein},
  volume~50 of {\em Princeton Math. Ser.}, pages 189--214. Princeton Univ.
  Press, Princeton, NJ, 2014.

\bibitem[Kes87a]{kesten-dla-hit-pro}
Harry Kesten.
\newblock Hitting probabilities of random walks on {${\bf Z}^d$}.
\newblock {\em Stochastic Process. Appl.}, 25(2):165--184, 1987.

\bibitem[Kes87b]{Kesten-phys}
Harry Kesten.
\newblock How long are the arms in {DLA}?
\newblock {\em J. Phys. A}, 20(1):L29--L33, 1987.

\bibitem[Kes90]{kesten-dla}
Harry Kesten.
\newblock Upper bounds for the growth rate of {DLA}.
\newblock {\em Phys. A}, 168(1):529--535, 1990.

\bibitem[Kig01]{kigami-book}
Jun Kigami.
\newblock {\em Analysis on fractals}, volume 143 of {\em Cambridge Tracts in
  Mathematics}.
\newblock Cambridge University Press, Cambridge, 2001.

\bibitem[KP04]{peres_krishnapur_collide}
M.~Krishnapur and Y.~Peres.
\newblock Recurrent graphs where two independent random walks collide finitely
  often.
\newblock {\em Electron. Commun. Probab.}, 9:72--81, 2004.

\bibitem[Law95]{lawler_1995}
G.~Lawler.
\newblock Subdiffusive fluctuations for internal diffusion limited aggregation.
\newblock {\em Ann. Probab.}, 23:71--86, 1995.

\bibitem[LBG92]{lawler_bramson_griffeath}
Gregory~F. Lawler, Maury Bramson, and David Griffeath.
\newblock Internal diffusion limited aggregation.
\newblock {\em Ann. Probab.}, 20(4):2117--2140, 1992.

\bibitem[Lev09]{levine-sandpile-tree}
Lionel Levine.
\newblock The sandpile group of a tree.
\newblock {\em European J. Combin.}, 30(4):1026--1035, 2009.

\bibitem[LP09]{peres_levine_strong_spherical}
Lionel Levine and Yuval Peres.
\newblock Strong spherical asymptotics for rotor-router aggregation and the
  divisible sandpile.
\newblock {\em Potential Anal.}, 30(1):1--27, 2009.

\bibitem[LP10]{peres_levine_scaling_limits}
Lionel Levine and Yuval Peres.
\newblock {S}caling {L}imits for {I}nternal {A}ggregation {M}odels with
  {M}ultiple {S}ources.
\newblock {\em Journal d'Analyse Math{\'e}matique}, 1:151--219, 2010.

\bibitem[LS18]{levine_silvestri-cylinder}
Lionel Levine and Vittoria Silvestri.
\newblock How long does it take for internal dla to forget its initial profile?
\newblock {\em Probability Theory and Related Fields}, 2018,
  \href{http://arxiv.org/abs/1801.08533}{{\ttfamily arXiv:1801.08533}}.
\newblock to appear.

\bibitem[Luc14]{lucas-drift-idla}
Cyrille Lucas.
\newblock The limiting shape for drifted internal diffusion limited aggregation
  is a true heat ball.
\newblock {\em Probab. Theory Related Fields}, 159(1-2):197--235, 2014.

\bibitem[Mar17]{martineau-directed-dla}
S\'{e}bastien Martineau.
\newblock Directed diffusion-limited aggregation.
\newblock {\em ALEA Lat. Am. J. Probab. Math. Stat.}, 14(1):249--270, 2017.

\bibitem[MD86]{meakin-deutch}
Paul Meakin and J.M. Deutch.
\newblock The formation of surfaces by diffusion limited annihilation.
\newblock {\em J. Chem. Phys.}, 85, 1986.

\bibitem[Osa90]{osada-pre-sierpinski}
Hirofumi Osada.
\newblock Isoperimetric constants and estimates of heat kernels of pre
  {S}ierpi\'{n}ski carpets.
\newblock {\em Probab. Theory Related Fields}, 86(4):469--490, 1990.

\bibitem[PRZ18]{dla-wedge}
Eviatar Procaccia, Ron Rosenthal, and Yuan Zhang.
\newblock Stabilization of {DLA} in a wedge.
\newblock 2018,  \href{http://arxiv.org/abs/1804.04236}{{\ttfamily
  arXiv:1804.04236}}.

\bibitem[PYZ19]{procaccia-scaling-limit}
Eviatar Procaccia, Jiayan Ye, and Yuan Zhang.
\newblock Stationary harmonic measure as the scaling limit of truncated
  harmonic measure.
\newblock 2019,  \href{http://arxiv.org/abs/1811.04793}{{\ttfamily
  arXiv:1811.04793}}.

\bibitem[PZ18]{procaccia-zero-harmonic}
Eviatar Procaccia and Yuan Zhang.
\newblock On sets of zero stationary harmonic measure.
\newblock 2018,  \href{http://arxiv.org/abs/1711.01013}{{\ttfamily
  arXiv:1711.01013}}.

\bibitem[PZ19]{Procaccia2019}
Eviatar~B. Procaccia and Yuan Zhang.
\newblock Stationary harmonic measure and {DLA} in the upper half plane.
\newblock {\em Journal of Statistical Physics}, pages 1572--9613, Jun 2019.

\bibitem[She10]{shellef_idla_percolation}
Eric Shellef.
\newblock I{DLA} on the supercritical percolation cluster.
\newblock {\em Electron. J. Probab.}, 15:no. 24, 723--740, 2010.

\bibitem[Spi76]{spitzer-principles}
Frank Spitzer.
\newblock {\em Principles of random walk}.
\newblock Springer-Verlag, New York-Heidelberg, second edition, 1976.
\newblock Graduate Texts in Mathematics, Vol. 34.

\bibitem[WS83]{dla}
T.~A. Witten and L.~M. Sander.
\newblock Diffusion-limited aggregation.
\newblock {\em Phys. Rev. B}, 27:5686--5697, May 1983.

\end{thebibliography}
\bibliographystyle{alpha_arxiv}
\end{document}